\documentclass[11pt]{article}

\usepackage{amsfonts}
\usepackage{stmaryrd}
\usepackage{latexsym}
\usepackage{amssymb}
\usepackage{mathrsfs}
\usepackage{dsfont}
\usepackage{amsmath}
\usepackage{amsthm}
\usepackage{cases}
\usepackage{extarrows}
\usepackage{algorithm}
\usepackage{algpseudocode}

\usepackage[german,english]{babel}

\usepackage[square, comma, sort&compress, numbers]{natbib}

\usepackage{fancyhdr}
\usepackage{geometry}
\usepackage{color}
\usepackage{epstopdf}
\usepackage{graphicx}
\usepackage{subfigure}

\usepackage{enumerate}

\theoremstyle{plain}
\newtheorem{theorem}{Theorem}[section]

\newtheorem{corollary}[theorem]{Corollary}
\newtheorem{lemma}[theorem]{Lemma}

\newtheorem{remark}[theorem]{Remark}
\newtheorem{definition}[theorem]{Definition}

\numberwithin{equation}{section}

\usepackage{hyperref}

\begin{document}

\title{Double determinants of quaternion matrices in the \(\mathrm{P}\)-generalized determinant framework with applications}

\author{{\ Jiale Gao$^{a}$,    Qing-Wen Wang$^{a,b}$\footnote{Corresponding author:  Qing-Wen Wang; Email: wqw@shu.edu.cn}}
\\ {\small$^{a}$Department of Mathematics, Shanghai University, Shanghai 200444, China}
\\
{\small $^b$National Key Laboratory of Maritime Unmanned Swarm Intelligence,} \\ {\small Shanghai University, Shanghai 200444, China}
}

\date{}

\maketitle

\begin{abstract}
Building on the \(\rm P\)-generalized determinant over quaternions,
we introduce the double determinant of arbitrary square quaternion matrices.
We then apply the double determinant to obtain a determinant representation for the inverse of any invertible quaternion matrix, without the Hermitian restriction.
Furthermore, we establish Cramer's rule for the unique solutions to the quaternion linear equations \(Ax=b\) and \(xA=b\), as well as to the quaternion matrix equations \(AX=B\), \(XA=B\), and \(AXB=C\).
As an application, we use the determinant representation of the unique solution to \(AXB=C\) to design a color image encryption--decryption scheme.
Numerical experiments demonstrate that the encryption algorithm effectively conceals the contour and color information of the original image, while the decryption algorithm restores the original image almost perfectly, thereby
confirming the effectiveness and practical potential of the proposed scheme.

\end{abstract}



\section{Introduction}
Quaternions, discovered by Hamilton in 1843 \cite{hamilton1853lectures,hamilton1899elements},  have been
wide practical applications in
color image processing \cite{grigoryan2018quaternion}, computer graphics \cite{goldman2010rethinking,vince2011quaternions}, scientific visualization \cite{hanson2006visualizing,hanson2024visualizing}, and aerospace engineering \cite{kuipers1999quaternions}.
Owing to the noncommutativity of quaternion multiplication, the classical determinant theory for complex matrices cannot be directly extended to quaternion matrices.
Consequently, how to define quaternion determinants and characterize their properties has long been a fundamental problem in quaternion matrix theory.

As early as 1845, Cayley \cite{cayley1845quat} first defined the determinant of a square quaternion matrix by expansion along the first column together with induction. However, the properties of the Cayley determinant are not satisfactory \cite{aslaksen1996quaternionic}. In 1920, Study \cite{study1920zur} established
 the quaternion determinant of a quaternion matrix via the determinant of its complex representation matrix, denoted by $\rm Sdet$.
In 1943, Dieudonn\'{e} \cite{dieudonne1943determinants} redefined the quaternion determinant by means of commutators, denoted by $\rm Ddet$. Furthermore, ${\rm Ddet}(A) = ({\rm Sdet}(A))^2$ for any  square quaternion  matrix $A$.
Later, Aslaksen \cite{aslaksen1996quaternionic} used the real representation matrix of a quaternion matrix to propose another quaternion determinant, following an approach similar to that of the Study determinant.
The Study, Dieudonn\'{e}, and Aslaksen determinants all satisfy the three axioms proposed by Aslaksen \cite{aslaksen1996quaternionic} for a quaternion determinant: the singularity axiom, the multiplicativity axiom, and the addition invariance axiom.
 However, Cohen and De Leo \cite{cohen2000quaternionic} pointed out that the Study determinant does not truly generalize the classical determinant of complex matrices: the former takes only non-negative real values, whereas the latter may take any complex value.  This limitation is also shared by the Dieudonn\'{e} and  Aslaksen determinants.

The pursuit of quaternion determinants taking values in all of $\mathbb{H}$ dates back to the early 1920s.
In 1922, Moore \cite{moore1922determinant} defined the determinant of quaternion Hermitian matrices by means of quaternion multivariate polynomials, denoted by $\rm Mdet$.
In 1980, Xie \cite{xie1980centralizable,xie1980expansion} gave another definition of the determinant for quaternion Hermitian matrices, denoted by $\rm Xdet$.
In 1991, Chen \cite{chen1991definition} proposed the quaternion determinant to arbitrary  square quaternion  matrices, denoted by $\rm Cdet$.
In 2008, Kyrchei \cite{kyrchei2008cramer,kyrchei2012theory} established the column-row determinant theory for quaternions,
in which the quaternion determinant can be expanded along any row or column, analogous to the classical Laplace expansion for complex matrix determinants; for details, see his recent monograph \cite{Kyrchei2026}.
In 2026, Gao and Wang \cite{PgeneralizedDetQuat} further introduced the $\rm P$-generalized determinant,  which  provides a unified framework for the Moore, Xie, and Chen determinants as well as the column-row determinant.
It not only significantly relaxes the stringent restrictions on the order of products imposed by these classical definitions but also retains most of their desirable properties.

However, the $\mathrm{P}$-generalized determinant introduced in \cite{PgeneralizedDetQuat} has thus far been applied only to quaternion Hermitian matrices, both for matrix inversion and for solving linear and matrix equations.
Inspired by Chen \cite{chen1991definition} and Kyrchei \cite{kyrchei2008cramer}, who respectively defined double determinants in terms of the Chen determinant and the column-row determinant, we introduce a double determinant for arbitrary square quaternion matrices within the
$\mathrm{P}$-generalized determinant framework.
The main contributions of this paper are summarized as follows:
\begin{enumerate}[(1)]
  \item By means of the double determinant and the \({\rm P}\)-generalized determinant,
        we derive a determinant representation for the inverse of an arbitrary invertible quaternion matrix.

  \item
  Using this determinant representation for invertible matrices,
we establish Cramer's rule for the unique solutions to the two classes of quaternion linear equations
\begin{equation*}
  Ax=b \quad \text{and} \quad xA=b,
\end{equation*}
as well as to the three classes of quaternion matrix equations
\begin{equation*}
  AX=B, \quad XA=B, \quad \text{and} \quad AXB=C,
\end{equation*}
where
the coefficient matrix is an arbitrary, not necessarily Hermitian, invertible quaternion matrix.

\item By combining Cramer's rule for \(AXB=C\), a block encryption--decryption strategy,
        and the Arnold transform together with its inverse,
        we design a color image encryption--decryption scheme.
        Numerical experiments demonstrate strong encryption performance
        and near-perfect recovery of the original image.

\end{enumerate}

The rest of this paper is organized as follows. Section~\ref{Preliminariessection} recalls the necessary notations, definitions, and lemmas. Section~\ref{Doubledetsection} defines the double determinant based on the \(\mathrm{P}\)-generalized determinant and provides determinant representations for matrix inverses.
Section~\ref{APPsection} establishes Cramer's rule for the quaternion  equations \(Ax=b\), \(xA=b\), \(AX=B\), \(XA=B\), and \(AXB=C\).
Section~\ref{numexpsection} describes a color image encryption--decryption scheme and presents numerical experiments verifying its effectiveness.
Section~\ref{conclussection} concludes the paper.

\section{Preliminaries}\label{Preliminariessection}

This section introduces the notations, definitions, and lemmas used throughout the paper. Let $\mathcal{Q}([n])=\{Q\subseteq[n]\mid Q\neq\emptyset\}$ be the set of all nonempty subsets of $[n]:=\{1,2,\ldots,n\}$, and let $S_n$ be the symmetric group on $[n]$.
For a permutation \(\sigma\in S_n\) with disjoint cycle decomposition \(\sigma=C_1C_2\cdots C_r\),
we define the sign of \(\sigma\) by
$
\operatorname{sgn}(\sigma)=(-1)^{n-r},
$
and denote the set of its cycles by \(\operatorname{cyc}(\sigma)=\{C_1,C_2,\dots,C_r\}\).
For a cycle $C\in S_n$, let $\operatorname{supp}(C)=\{i\in[n]\mid C(i)\neq i\}\subseteq[n]$ be the support of $C$,
and let $|\operatorname{supp}(C)|$ be the length of $C$. Let the quaternion algebra be
\begin{equation*}
\mathbb{H}=\{q_1+q_2\mathbf{i}+q_3\mathbf{j}+q_4\mathbf{k}\mid \mathbf{i}^2=\mathbf{j}^2=\mathbf{k}^2=\mathbf{ijk}=-1,\ q_1,q_2,q_3,q_4\in\mathbb{R}\},
\end{equation*}
where $\mathbb{R}$ denotes the field of real numbers.
For a quaternion $q=q_0+q_1\mathbf{i}+q_2\mathbf{j}+q_3\mathbf{k}$ with $q_0,q_1,q_2,q_3\in\mathbb{R}$, let $\overline{q}=q_0-q_1\mathbf{i}-q_2\mathbf{j}-q_3\mathbf{k}$ be the conjugate of $q$, and $|q|=\sqrt{q\overline{q}}=\sqrt{q_0^2+q_1^2+q_2^2+q_3^2}$ be the norm of $q$. Let $\mathbb{H}^{m\times n}$ denote the set of all $m\times n$ quaternion matrices.
For $A\in\mathbb{H}^{m\times n}$, let $\overline{A}$, $A^{T}$, and $A^{*}$ denote the conjugate, transpose, and conjugate transpose of $A$, respectively;
let $\operatorname{rank}_{r\text{-}\mathrm{col}}(A) $, ${\rm rank}_{l\text{-}{\rm col}}(A)$, ${\rm rank}_{r\text{-}{\rm row}}(A)$, and ${\rm rank}_{l\text{-}{\rm row}}(A)$ denote the right column rank, left  column rank,  right row  rank, and left row   rank of $A$, respectively;
 let $A^{\alpha}_{\beta} \in\mathbb{H}^{m_1\times n_1}$ denote the submatrix of $A$ formed by the rows indexed by $\alpha=\{\alpha_1,\ldots,\alpha_{m_1}\}\subseteq[m]$ and the columns indexed by $\beta=\{\beta_1,\ldots,\beta_{n_1}\}\subseteq[n]$;
  let $A^{ij}\in\mathbb{H}^{(m-1)\times(n-1)}$ denote the matrix obtained from $A$ by deleting the $i$-th row and the $j$-th column;
   let $A_{\cdot j}(b)\in\mathbb{H}^{m\times n}$ denote the matrix obtained by replacing the $j$-th column with the column vector $b\in\mathbb{H}^{m\times 1}$;
    let $A_{i\cdot}(c)\in\mathbb{H}^{m\times n}$ denote the matrix obtained from $A$ by replacing the $i$-th row with the row vector $c\in\mathbb{H}^{1\times n}$;
     let $a_{i\cdot}\in\mathbb{H}^{1\times n}$, $a_{\cdot j}\in\mathbb{H}^{m\times 1}$, $a_{i\cdot}^{*}\in\mathbb{H}^{1\times m}$, and $a_{\cdot j}^{*}\in\mathbb{H}^{n\times 1}$ denote the $i$-th row of $A$, the $j$-th column of $A$, the $i$-th row of $A^{*}$, and the $j$-th column of $A^{*}$, respectively.
  For any quaternion matrix $B\in\mathbb{H}^{m\times n}$, let $B_{i\cdot}\in\mathbb{H}^{1\times n}$ and $B_{\cdot j}\in\mathbb{H}^{m\times 1}$ denote the $i$-th row and the $j$-th column of $B$, respectively. Let $0$ denote the zero matrix of appropriate size and $I_n$ the $n\times n$ identity matrix.
  Let $P_{ij}(b)$ denote the elementary addition matrix, obtained from $I_n$ by adding a left multiple $b\in\mathbb{H}$ of the $j$-th row to the $i$-th row;
  let $P_i(b)$ denote the elementary scaling matrix, obtained from $I_n$ by left-multiplying the $i$-th row by $b\in\mathbb{H}$;
  and let $P_{ij}$ denote the elementary interchange matrix, obtained from $I_n$ by interchanging the $i$-th and $j$-th rows.
   Throughout this paper, we adopt the following convention: for any matrix and its submatrices, the row and column indices are taken from the same index set; that is, a submatrix inherits the row and column indices of the original matrix and is not reindexed.

\begin{definition}\cite{PgeneralizedDetQuat}
Let $S_n$ be the symmetric group on the set $[n]$.
\begin{enumerate}[(1)]
  \item
Define a function
\begin{equation*}
  {\mathrm{sel}}:\mathcal{Q}([n])\to[n],
\end{equation*}
such that ${\mathrm{sel}}(Q)\in Q$ for every nonempty $Q\in\mathcal{Q}([n])$.
We call this function a selector on $\mathcal{Q}([n])$ and refer to ${\mathrm{sel}}(Q)$ as the selected element of $Q$.

\item
Let $\vartriangle$ be a strict total order on $[n]$, and let $  {\mathrm{sel}}$ be a selector on $\mathcal{Q}([n])$.  For
  $\sigma\in S_n$,
define a strict total order $\mathrel{\curlywedge}_{\vartriangle}$  on $\mathrm{cyc}(\sigma)$ by
\begin{equation*}
  C\,\mathrel{\curlywedge}_{\vartriangle}\,D \iff
  {\rm sel}(\mathrm{supp}(C))\,  \vartriangle  \,{\rm sel}(\mathrm{supp}(D)),
\end{equation*}
where $C,D\in\mathrm{cyc}(\sigma)$ are two distinct cycles.
 To emphasize its dependence on $\mathrm{cyc}(\sigma)$, $\vartriangle$ and $  {\mathrm{sel}}$, we also write this order as ${\curlywedge}_{\vartriangle}\big(\mathrm{cyc}(\sigma), \vartriangle,  {\mathrm{sel}}\big)$.

\end{enumerate}
\end{definition}

\begin{definition}\cite{PgeneralizedDetQuat}
Let $S_n$ be the symmetric group on the set $[n]$.
\begin{enumerate}[(1)]
\item
 Let $ {\mathrm{first}}$ be a selector on  $\mathcal{Q}([n])$,
 let  $\triangleleft$ be a strict total order on $[n]$, and let
$\prec_{\triangleleft}:= \prec_{\triangleleft}\big(\mathrm{cyc}(\sigma), \triangleleft, {\mathrm{first}}  \big) $ for $\sigma\in S_n$.
 Define  the left-order representation of the cycle $C\in\mathrm{cyc}(\sigma)$  based on the selector  ${\rm first}$ by
\begin{equation*}
 C_{\mathfrak{l}}= \big(c,
\sigma(c), \sigma^2(c), \dots ,
\sigma^{|\mathrm{supp}(C)|-1}(c)\big),
\end{equation*}
where $c={\mathrm{first}}(\mathrm{supp}(C))$.
 Furthermore, define the left-order cycle representation of the permutation $\sigma$ based on the selector ${\rm first}$ and the strict total order $\prec_{{\triangleleft}} $ by
\begin{equation}\label{cicdecsigmapereq}
  \sigma_l={(C_1)}_{\mathfrak{l}} (C_2)_{\mathfrak{l}}
\dots(C_k)_{\mathfrak{l}},
\end{equation}
where $C_1,\dots,C_k$ are all cycles of $\sigma$, arranged in increasing order
with respect to $\prec_{\triangleleft}$, i.e.,
\begin{equation*}
C_1\prec_{\triangleleft} C_2\prec_{\triangleleft}\cdots\prec_{\triangleleft} C_k.
\end{equation*}
 The rule that determines the left-order cycle representation  \eqref{cicdecsigmapereq}  of every permutation in $S_n$, based on the selector ${\rm first}$ and the strict total order $\prec_{{\triangleleft}}$, is called a left-order cycle representation rule. It is denoted by ${\rm P}_l(S_n,{\rm first},\prec_{{\triangleleft}})$ and abbreviated as ${\rm P}_l$.
    In particular, if, for every $\sigma\in S_n$, the selected element of the cycle $(C_1)_\mathfrak{l}$ in \eqref{cicdecsigmapereq} is fixed at a prescribed element $i\in[n]$, then this rule is denoted by ${\rm P}_{l\text{-}i}$.

\item
 Let $ {\mathrm{last}}$ be a selector on  $\mathcal{Q}([n])$,
let  $\triangleright$  be a  strict total order   on $[n]$,
 and let $\succ_{{\triangleright}}:=
 \succ_{{\triangleright}}
 \big(\mathrm{cyc}(\sigma),
 \triangleright, {\mathrm{last}}  \big) $ for $\sigma\in S_n$.
Define the right-order representation of the cycle $D \in\mathrm{cyc}(\sigma)$ based on the selector  ${\rm last}$ by
\begin{equation*}
 D_{\mathfrak{r}}= \big(
\sigma^{|\mathrm{supp}(D)|-1}(d) ,\dots , \sigma^2(d) , \sigma(d) , d  \big),
\end{equation*}
where $d={\mathrm{last}}(\mathrm{supp}(D))$.
 Furthermore, define the right-order cycle representation of the permutation $\sigma$ based on the selector ${\rm last}$ and the strict total order $\succ_{{\triangleright}}$ by
\begin{equation}\label{sigmadecrighteq}
  \sigma_r=(D_k)_{\mathfrak{r}} \dots (D_2)_{\mathfrak{r}} {(D_1)}_{\mathfrak{r}} ,
\end{equation}
where $D_1,\dots,D_k $ are all cycles of $\sigma$,
 arranged in decreasing order with respect to $\succ_{{\triangleright}}$, i.e.,
 \begin{equation*}
 D_k\succ_{{\triangleright}}  \ldots \succ_{{\triangleright}}
 D_2\succ_{{\triangleright}} D_1.
\end{equation*}
The rule that determines the right-order cycle representation \eqref{sigmadecrighteq} of every permutation in $S_n$, based on the selector ${\rm last}$ and the strict total order $\succ_{{\triangleright}}$, is called a right-order cycle representation rule. It is denoted by ${\rm P}_r(S_n,{\rm last},\succ_{{\triangleright}})$ and abbreviated as ${\rm P}_r$.
In particular, if, for every $\sigma\in S_n$, the selected element of the cycle $(D_1)_\mathfrak{r}$ in \eqref{sigmadecrighteq} is fixed at a prescribed element $j\in[n]$, then this  rule is denoted by ${\rm P}_{r\text{-}j}$.

\item
The left-order cycle representation rule ${\rm P}_l$ and the right-order cycle representation rule ${\rm P}_r$ are collectively referred to as cycle representation rules (or simply rules) and denoted generically by ${\rm P}$.

\end{enumerate}

\end{definition}

Unless otherwise stated, in this paper $\rm P$, ${\rm P}_l$, and ${\rm P}_r$  denote the cycle representation rule, the  left-order cycle representation rule, and the right-order cycle representation rule, respectively.

\begin{definition}\label{Pgdetdeinition}\cite{PgeneralizedDetQuat}
Let $S_n$ be the symmetric group on  $[n]$, and  $A=[a_{ij}]\in\mathbb{H}^{n\times n}$.
\begin{enumerate}[(1)]
  \item
   Let ${\rm P}_l={\rm P}_l( S_n, {\rm first}, \prec_{{\triangleleft}})$ be a left-order cycle representation rule.
Under the rule ${\rm P}_l$, suppose that the left-order cycle representation  of a permutation $\sigma\in S_n$ is
\begin{equation}\label{sigmmaldeceq}
  \sigma_l=(i_{k_1},i_{k_1+1},...,i_{k_1+v_1})(i_{k_2},i_{k_2+1},...,i_{k_2+v_2})...
  (i_{k_u},i_{k_u+1},...,i_{k_u+v_u}),
\end{equation}
where $i_{k_1},...,i_{k_1+v_1},...,i_{k_u},...,i_{k_u+v_u}\in[n]$ are fixed and mutually distinct. Let
\begin{equation}\label{fasiggmaldeceq}
  f_A(\sigma_l)=
  (a_{i_{k_1},i_{k_1+1}}a_{i_{k_1+1},i_{k_1+2}}...a_{i_{k_1+v_1},i_{k_1}})
\cdots
  (a_{i_{k_u},i_{k_u+1}}a_{i_{k_u+1},i_{k_u+2}}...a_{i_{k_u+v_u},i_{k_u}}).
\end{equation}
Then the ${\rm P}_l$-generalized determinant of $A$, denoted by
$\mathrm{GDet}_{{\rm P}_l}(A)$, is defined by
\begin{equation*}
  \mathrm{GDet}_{{\rm P}_l}(A)=
  \sum_{\sigma\in S_n}\mathrm{sgn}(\sigma)f_A(\sigma_l).
\end{equation*}

  \item
  Let ${\rm P}_r={\rm P}_r( S_n, {\rm last}, \succ_{{\triangleright}})$ be a right-order cycle representation rule.
      Under the rule ${\rm P}_{r}$, suppose that the right-order cycle representation  of a permutation $\sigma\in S_n$ is
\begin{equation*}
  \sigma_r=(j_{k_u+v_u},...,j_{k_u+1},j_{k_u})...(j_{k_2+v_2},...,j_{k_2+1},j_{k_2})(j_{k_1+v_1},...,j_{k_1+1},j_{k_1}),
\end{equation*}
where $j_{k_u+v_u},...,j_{k_u},...,j_{k_1+v_1},...,j_{k_1}\in[n]$ are fixed and mutually distinct.
Let
\begin{equation*}
  f_A( \sigma_r)= (a_{j_{k_u},j_{k_u+ v_u}}...a_{j_{k_u+2},j_{k_u+1}}a_{j_{k_u+1},j_{k_u}})\cdots
(a_{j_{k_1},j_{k_1+ v_1}}...a_{j_{k_1+2},j_{k_1+1}}a_{j_{k_1+1},j_{k_1}}).
\end{equation*}
Then the ${\rm P}_{r}$-generalized determinant of $A$,
denoted by
$ {\rm GDet}_{{\rm P}_r} (A)$,
is defined by
\begin{equation*}
 {\rm GDet}_{{\rm P}_r}
 (A)=\sum\limits_{\sigma\in S_n} {\rm sgn}(\sigma)   f_A( \sigma_r).
\end{equation*}

\item The ${\rm P}_l$-generalized determinant ${\rm GDet}_{{\rm P}_l}(A)$ and the ${\rm P}_r$-generalized determinant ${\rm GDet}_{{\rm P}_r}(A)$  are collectively called the ${\rm P}$-generalized determinants, denoted by ${\rm GDet}_{{\rm P}}(A)$.

\end{enumerate}

\end{definition}

\begin{definition}\label{varimdetremark}
\begin{enumerate}[(1)]
\item\cite{moore1922determinant}   If ${\rm first}(Q)=\min(Q)$ for $Q\in \mathcal{Q}([n])$, and  the strict total order $\triangleleft$  on $[n]$  is defined by
$
n\triangleleft n-1\triangleleft \ldots \triangleleft 2\triangleleft 1,
$
then $ {\rm Mdet}(A):={\rm GDet}_{{\rm P}_l}(A) $  for a Hermitian matrix $A\in\mathbb{H}^{n\times n}$.

\item\cite{xie1980expansion} If ${\rm first}(Q)=\min(Q)$ for $Q\in \mathcal{Q}([n])$, and the strict total order $\triangleleft$ on $[n]$ is  defined by
$1\triangleleft 2\triangleleft \ldots \triangleleft n$, then
${\rm Xdet}(A):={\rm GDet}_{{\rm P}_l}(A)$ for a Hermitian matrix $A\in\mathbb{H}^{n\times n}$.

\item\cite{chen1991definition}
If ${\rm first}(Q)=\max(Q)$ for $Q\in \mathcal{Q}([n])$, and the strict total order $\triangleleft$ on $[n]$ is defined by
$n \triangleleft n-1\triangleleft \ldots \triangleleft 2\triangleleft 1$,
then ${\rm Cdet}(A):={\rm GDet}_{{\rm P}_l}(A)$ for $A\in\mathbb{H}^{n\times n}$.

\item\cite{kyrchei2008cramer}\label{PGdetspecaseremarkitem4}
Let $i\in[n]$.
If
\begin{equation*}
  {\rm first}(Q)= \begin{cases}
i, & i\in Q,\\
\min(Q), & i\notin Q
\end{cases}
\end{equation*}
for $Q\in \mathcal{Q}([n])$, and the strict total order $\triangleleft$ on $[n]$ is defined by
$i\triangleleft i_1 \triangleleft i_2 \triangleleft \ldots \triangleleft i_{n-1}$,
where $i_1, \ldots, i_{n-1} \in [n]\setminus \{i\}$ satisfy
$i_1 < i_2 < \ldots < i_{n-1}$,
then ${\rm rdet}_i(A):={\rm GDet}_{{\rm P}_{l\text{-}i}}(A)$ for $A\in\mathbb{H}^{n\times n}$.

\item\cite{kyrchei2008cramer}\label{PGdetspecaseremarkitem5}
Let $j\in[n]$.  If
\begin{equation*}
  {\rm last}(Q)= \begin{cases}
j, & j\in Q,\\
\min(Q), & j\notin Q
\end{cases}
\end{equation*}
for $Q\in \mathcal{Q}([n])$, and the strict total order $\triangleright$ on $[n]$ is defined by
$j_{n-1}\triangleright \ldots \triangleright j_2 \triangleright j_1 \triangleright j$,
where $j_{n-1}, \ldots, j_1 \in [n]\setminus \{j\}$ satisfy
$j_{n-1}>\ldots>j_2>j_1$,
then ${\rm cdet}_j(A):={\rm GDet}_{{\rm P}_{r\text{-}j}}(A)$ for  $A\in\mathbb{H}^{n\times n}$.

\item\cite{kyrchei2008cramer}
If $A\in\mathbb{H}^{n\times n}$ is Hermitian, then we denote  ${\rm det}(A):={\rm cdet}_j(A)={\rm rdet}_i(A)$ for any $i,j\in[n]$.
\end{enumerate}

\end{definition}

\begin{lemma}\cite[Theorem 3.9]{PgeneralizedDetQuat}\label{basicproPgdecThle}
Let $A=[a_{xy}]\in\mathbb{H}^{n\times n}$ for all $x,y\in[n]$, and let $q \in \mathbb{H}$.
\begin{enumerate}[(1)]
  \item
$
     {\rm GDet}_{{\rm P}_{l\text{-}i}}\left(A_{i\cdot}(q\cdot a_{i\cdot})\right)
    = q\cdot {\rm GDet}_{{\rm P}_{l\text{-}i}}(A), ~
      {\rm GDet}_{{\rm P}_{r\text{-}j}}\left(A_{\cdot j}(a_{\cdot j}\cdot q)\right)
    = {\rm GDet}_{{\rm P}_{r\text{-}j}}(A)\cdot q.
$

  \item\label{basicproPgdecThitem3}
If $t\in [n]$  and  $a_{t \cdot }=b+c$ with
    $b,c \in \mathbb{H}^{1 \times n}$,     then
  \[
    {\rm GDet}_{{\rm P}}(A)
    = {\rm GDet}_{{\rm P}}(A_{t\cdot}(b))
    + {\rm GDet}_{{\rm P}}(A_{t\cdot}(c)).
  \]

  \item
 If $t\in [n]$  and $a_{\cdot t}=b+c$
  with $b,c\in \mathbb{H}^{n \times 1}$, then
  \[
    {\rm GDet}_{{\rm P}}(A)
    = {\rm GDet}_{{\rm P}}(A_{\cdot t}(b))
    + {\rm GDet}_{{\rm P}}(A_{\cdot t}(c)).
  \]
\end{enumerate}
\end{lemma}

\begin{lemma}\cite[Theorem 4.1]{PgeneralizedDetQuat}\label{HerPGdetmainprothle}
Let $A\in\mathbb{H}^{n\times n}$ be Hermitian. Then ${\rm GDet}_{{\rm P}}(A)$ is a real number  independent of the cycle representation rule  ${\rm P}$, denoted  by  ${\rm Det}(A)$.
\end{lemma}

\section{Double determinant}\label{Doubledetsection}

To study the P-generalized determinant for arbitrary quaternion matrices (not necessarily Hermitian), we first give the following properties:

\begin{theorem}\label{DetBgBAgAvarTh}
Let $A \in \mathbb{H}^{n\times n}$ and $b \in \mathbb{H}$. Then the following statements hold:
\begin{enumerate}[(1)]
  \item\label{DetBgBAgAvarThitem1} If $B=AP_{ij}(b)$, then
  $
  \operatorname{Det}(B^*B)=\operatorname{Det}(A^*A).
  $

  \item\label{DetBgBAgAvarThitem2} If $B=AP_{j}(b)$, then
  $
  \operatorname{Det}(B^*B)=|b|^{2}\operatorname{Det}(A^*A).
  $

  \item\label{DetBgBAgAvarThitem3} If $B=AP_{ij}$, then
  $
  \operatorname{Det}(B^*B)=\operatorname{Det}(A^*A).
  $

  \item\label{DetBgBAgAvarThitem4}
  If  \(b_k\in\mathbb{H}\)  for all \(k\in[n]\setminus\{j\}\) and
  $
  B=A_{\cdot j}\left(\sum_{k=1,\,k\neq j}^{n} a_{\cdot k} \cdot b_k \right)
  $,
 then
  $
  \operatorname{Det}(B^*B)=0.
  $

  \item\label{DetBgBAgAvarThitem5}

  If    $b_k \in \mathbb{H}$ for  all $k\in[n]$ and
  $
  B=A_{\cdot j}\left(\sum_{k=1}^{n} a_{\cdot k} \cdot b_k\right)
  $,
then
  $
  \operatorname{Det}(B^*B)=|b_j|^{2}\operatorname{Det}(A^*A).
  $
\end{enumerate}
\end{theorem}

\begin{proof}
Items \eqref{DetBgBAgAvarThitem1}--\eqref{DetBgBAgAvarThitem3} follow immediately from Items (1)--(3) in \cite[Remark 4.11]{PgeneralizedDetQuat}, while Item \eqref{DetBgBAgAvarThitem4} is a consequence of \cite[Corollary 4.9]{PgeneralizedDetQuat}.

It remains to prove Item \eqref{DetBgBAgAvarThitem5}. Indeed, let $A=[\alpha_1, \alpha_2,\ldots, \alpha_n]$ and $\beta_j = \sum_{k=1}^{n}\alpha_k b_k$. Then
\begin{equation*}
  B= A_{\cdot j}\left(\sum_{k=1}^{n} a_{\cdot k} \cdot b_k \right)=[\alpha_1, \ldots,  \alpha_{j-1}, \beta_j, \alpha_{j+1}, \ldots,  \alpha_n   ].
\end{equation*}
Let $B^*B=[\gamma_1, \gamma_2,\ldots, \gamma_n]$ and $A^*A= [ \delta_1^{T}, \delta_2^{T},\ldots, \delta_n^{T}] ^{T}$.
Using  Lemma \ref{basicproPgdecThle} and \cite[Corollary 4.9]{PgeneralizedDetQuat}, we have
\begin{align*}
{\rm GDet}_{{\rm P}_{r\text{-}j}}
 \left[
  \begin{array}{ccccc}
  \alpha_1^* \alpha_1   & \ldots &  \alpha_1^* (\sum _{k=1,k\neq j}^{n} \alpha_k  b_k ) & \ldots &  \alpha_1^* \alpha_n \\
    \vdots  & \ldots &  \vdots & \ldots &  \vdots \\
  \beta_j^* \alpha_1   & \ldots &  \beta_j^* (\sum _{k=1,k\neq j}^{n} \alpha_k  b_k )& \ldots &  \beta_j^* \alpha_n \\
    \vdots  & \ldots &  \vdots & \ldots &  \vdots \\
  \alpha_n^* \alpha_1   & \ldots & \alpha_n^*(\sum _{k=1,k\neq j}^{n} \alpha_k  b_k )& \ldots &  \alpha_n^* \alpha_n \\
  \end{array}
\right]
={\rm GDet}_{{\rm P}_{r\text{-}j}} (B^*B)_{\cdot j}\left( \sum_{k=1,k\neq j}^{n} \gamma_k b_k \right)=0,
\end{align*}
\begin{align*}
 &{\rm GDet}_{{\rm P}_{r\text{-}j}}\left[
  \begin{array}{ccccc}
  \alpha_1^* \alpha_1   & \ldots &  \alpha_1^*   \alpha_j  b_j  & \ldots &  \alpha_1^* \alpha_n \\
    \vdots  & \ldots &  \vdots & \ldots &  \vdots \\
  \beta_j^* \alpha_1   & \ldots &  \beta_j^* \alpha_j  b_j & \ldots &  \beta_j^* \alpha_n \\
    \vdots  & \ldots &  \vdots & \ldots &  \vdots \\
  \alpha_n^* \alpha_1   & \ldots &  \alpha_n^*  \alpha_j  b_j & \ldots &  \alpha_n^* \alpha_n \\
  \end{array}
\right]=
 \left({\rm GDet}_{{\rm P}_{r\text{-}j}}\left[
  \begin{array}{ccccc}
  \alpha_1^* \alpha_1   & \ldots &  \alpha_1^*   \alpha_j   & \ldots &  \alpha_1^* \alpha_n \\
    \vdots  & \ldots &  \vdots & \ldots &  \vdots \\
  \beta_j^* \alpha_1   & \ldots &  \beta_j^* \alpha_j   & \ldots &  \beta_j^* \alpha_n \\
    \vdots  & \ldots &  \vdots & \ldots &  \vdots \\
  \alpha_n^* \alpha_1   & \ldots &  \alpha_n^*  \alpha_j   & \ldots &  \alpha_n^* \alpha_n \\
  \end{array}
\right]\right) b_j\\
=& \left({\rm GDet}_{{\rm P}_{r\text{-}j}}\left[
  \begin{array}{ccccc}
  \alpha_1^* \alpha_1   & \ldots &  \alpha_1^*   \alpha_j   & \ldots &  \alpha_1^* \alpha_n \\
    \vdots  & \ldots &  \vdots & \ldots &  \vdots \\
  (\sum_{k=1}^{n} \overline{b_k} \alpha_k^*) \alpha_1   & \ldots &  (\sum_{k=1}^{n} \overline{b_k} \alpha_k^*) \alpha_j   & \ldots &  (\sum_{k=1}^{n} \overline{b_k} \alpha_k^*) \alpha_n \\
    \vdots  & \ldots &  \vdots & \ldots &  \vdots \\
  \alpha_n^* \alpha_1   & \ldots &  \alpha_n^*  \alpha_j   & \ldots &  \alpha_n^* \alpha_n \\
  \end{array}
\right]\right) b_j\\
=&  \left( {\rm GDet}_{{\rm P}_{r\text{-}j}}  (A^*A)_{j\cdot} \left( \sum_{k=1}^{n}  \overline{b_k}   \delta_k  \right) \right)  b_j
= \overline{b_j}  {\rm Det}(A^*A) b_j= |b_j|^{2}\operatorname{Det}(A^*A).
\end{align*}
Then, from Lemmas  \ref{HerPGdetmainprothle}  and \ref{basicproPgdecThle}, it follows
\begin{align*}
&{\rm Det}(B^*B) =
 {\rm Det} \left[
  \begin{array}{ccccc}
  \alpha_1^* \alpha_1   & \ldots &  \alpha_1^* \beta_j & \ldots &  \alpha_1^* \alpha_n \\
    \vdots  & \ldots &  \vdots & \ldots &  \vdots \\
  \beta_j^* \alpha_1   & \ldots &  \beta_j^* \beta_j & \ldots &  \beta_j^* \alpha_n \\
    \vdots  & \ldots &  \vdots & \ldots &  \vdots \\
  \alpha_n^* \alpha_1   & \ldots &  \alpha_n^* \beta_j & \ldots &  \alpha_n^* \alpha_n \\
  \end{array}
\right]
={\rm GDet}_{{\rm P}_{r\text{-}j}} \left[
  \begin{array}{ccccc}
  \alpha_1^* \alpha_1   & \ldots &  \alpha_1^* \beta_j & \ldots &  \alpha_1^* \alpha_n \\
    \vdots  & \ldots &  \vdots & \ldots &  \vdots \\
  \beta_j^* \alpha_1   & \ldots &  \beta_j^* \beta_j & \ldots &  \beta_j^* \alpha_n \\
    \vdots  & \ldots &  \vdots & \ldots &  \vdots \\
  \alpha_n^* \alpha_1   & \ldots &  \alpha_n^* \beta_j & \ldots &  \alpha_n^* \alpha_n \\
  \end{array}
\right]\\
=&
{\rm GDet}_{{\rm P}_{r\text{-}j}}
 \left[
  \begin{array}{ccccc}
  \alpha_1^* \alpha_1   & \ldots &  \alpha_1^* (\sum _{k=1,k\neq j}^{n} \alpha_k  b_k ) & \ldots &  \alpha_1^* \alpha_n \\
    \vdots  & \ldots &  \vdots & \ldots &  \vdots \\
  \beta_j^* \alpha_1   & \ldots &  \beta_j^* (\sum _{k=1,k\neq j}^{n} \alpha_k  b_k )& \ldots &  \beta_j^* \alpha_n \\
    \vdots  & \ldots &  \vdots & \ldots &  \vdots \\
  \alpha_n^* \alpha_1   & \ldots & \alpha_n^*(\sum _{k=1,k\neq j}^{n} \alpha_k  b_k )& \ldots &  \alpha_n^* \alpha_n \\
  \end{array}
\right] +{\rm GDet}_{{\rm P}_{r\text{-}j}}\left[
  \begin{array}{ccccc}
  \alpha_1^* \alpha_1   & \ldots &  \alpha_1^*   \alpha_j  b_j  & \ldots &  \alpha_1^* \alpha_n \\
    \vdots  & \ldots &  \vdots & \ldots &  \vdots \\
  \beta_j^* \alpha_1   & \ldots &  \beta_j^* \alpha_j  b_j & \ldots &  \beta_j^* \alpha_n \\
    \vdots  & \ldots &  \vdots & \ldots &  \vdots \\
  \alpha_n^* \alpha_1   & \ldots &  \alpha_n^*  \alpha_j  b_j & \ldots &  \alpha_n^* \alpha_n \\
  \end{array}
\right],
\end{align*}
which, together with  the above two equations, shows
$ {\rm Det}(B^*B)= |b_j|^{2}\operatorname{Det}(A^*A)$.
\end{proof}

\begin{remark}
In view of the invariance of the
 $\mathrm{P}$-generalized determinant for quaternion Hermitian matrices established in Lemma \ref{HerPGdetmainprothle},
Items \eqref{DetBgBAgAvarThitem2} and \eqref{DetBgBAgAvarThitem4} of Theorem \ref{DetBgBAgAvarTh} also follow directly from \cite[Theorems 8.1 and 8.7]{kyrchei2012theory}, respectively.
\end{remark}

In parallel with the column operations in Theorem~\ref{DetBgBAgAvarTh}, we now state the corresponding results for row operations.

\begin{theorem}\label{DetBBgAAgvarTh}
Let $A \in \mathbb{H}^{n\times n}$ and $b \in \mathbb{H}$. Then the following statements hold:
\begin{enumerate}[(1)]
  \item\label{DetBBgAAgvarThitem1} If $B=P_{ij}(b)A$, then
  $
  \operatorname{Det}(BB^*)=\operatorname{Det}(AA^*).
  $

  \item If $B=P_{j}(b)A$, then
  $
  \operatorname{Det}(BB^*)=|b|^{2}\operatorname{Det}(AA^*).
  $

  \item If $B=P_{ij}A$, then
  $
  \operatorname{Det}(BB^*)=\operatorname{Det}(AA^*).
  $

  \item If $b_k \in \mathbb{H}$ for all $k\in[n]\setminus\{i\}$ and
  $
  B=A_{ i \cdot}\left(\sum_{k=1,\,k\neq i}^{n}   b_k \cdot a_{k \cdot } \right)
  $, then
  $
  \operatorname{Det}(BB^*)=0.
  $

  \item  If $b_k \in \mathbb{H}$ for all $k\in[n]$ and
  $
  B=A_{ i \cdot}\left(\sum_{k=1}^{n}b_k  \cdot  a_{ k \cdot} \right)
  $, then
  $
  \operatorname{Det}(BB^*)=|b_i|^{2}\operatorname{Det}(AA^*).
  $
\end{enumerate}
\end{theorem}

\begin{proof}
The proof is analogous to that of Theorem \ref{DetBgBAgAvarTh}, and we omit the details.
\end{proof}

Since every principal submatrix of a Hermitian matrix is also Hermitian, Kyrchei \cite[Definition 9.5]{kyrchei2012theory} defined, via row-column determinants, the determinant rank of a Hermitian matrix \(A\), denoted by \(\mathrm{rank}_{d}(A)\); see also \cite[Definition 3.2.18]{Kyrchei2026}. Inspired by this, we next introduce the concept of the determinant rank for a quaternion Hermitian matrix with respect to the $\mathrm{P}$-generalized determinant,
and discuss its relations to the left row, right row, left column, and right column ranks.

\begin{definition}
Let $A\in\mathbb{H}^{n\times n}$ be Hermitian.
The determinant rank of $A$, denoted by ${\rm rank}_{D}(A)$, is defined as the largest integer $r$ such that $A$ has a nonzero principal minor of order $r$.
\end{definition}

\begin{corollary}\label{rankDdrowcolrcrankTH}
Let $A \in \mathbb{H}^{m \times n}$. Then
\begin{align}
&  {\rm rank}_{D}(A^*A) = {\rm rank}_{d}(A^*A)=  {\rm rank}_{r\text{-}{\rm col}}(A^*A)=
  {\rm rank}_{r\text{-}{\rm col}}(A)=  {\rm rank}_{l\text{-}{\rm row}}(A^*)
 \label{rankDdrowcolrcrankTHitem01}  \\
  =&{\rm rank}_{D}(AA^*)  = {\rm rank}_{d}(AA^*)=  {\rm rank}_{l\text{-}{\rm row}}(AA^*)
  =  {\rm rank}_{l\text{-}{\rm row}}(A)=  {\rm rank}_{r\text{-}{\rm col}}(A^*)  \label{rankDdrowcolrcrankTHitem02}  \\
=&  {\rm rank}_{l\text{-}{\rm col}}(A^T)=  {\rm rank}_{l\text{-}{\rm col}}(\overline{A})=
 {\rm rank}_{r\text{-}{\rm row}}(A^T)=  {\rm rank}_{r\text{-}{\rm row}}(\overline{A})   \label{rankDdrowcolrcrankTHitem03}.
\end{align}
\end{corollary}

\begin{proof}
Using  Lemma  \ref{HerPGdetmainprothle} and \cite[Theorems 9.7, 9.8]{kyrchei2012theory} (or \cite[Theorems 3.2.24, 3.2.25]{Kyrchei2026}), we obtain \eqref{rankDdrowcolrcrankTHitem01} and \eqref{rankDdrowcolrcrankTHitem02} directly.
By standard rank relations for quaternion matrices,
\begin{align*}
  & {\rm rank}_{l\text{-}{\rm row}}(A)    ={\rm rank}_{l\text{-}{\rm col}}(A^T)
  , \quad
  {\rm rank}_{r\text{-}{\rm col}}(A)    ={\rm rank}_{r\text{-}{\rm row}}(A^T),
\end{align*}
we obtain \eqref{rankDdrowcolrcrankTHitem03}.
\end{proof}

\begin{remark}
By Corollary~\ref{rankDdrowcolrcrankTH}, the quantities in \eqref{rankDdrowcolrcrankTHitem01}--\eqref{rankDdrowcolrcrankTHitem03} are equal. Their common value is called the rank of the quaternion matrix \(A\in\mathbb{H}^{m\times n}\), and is denoted by \(\operatorname{rank}(A)\).
\end{remark}

Since \(\operatorname{det}(AA^*)=\operatorname{det}(A^*A)\) for any \(A\in\mathbb{H}^{n\times n}\) by \cite[Theorem 7.4]{kyrchei2008cramer}, Kyrchei \cite[Definition 7.2]{kyrchei2008cramer} defined the double determinant by
\[
\operatorname{ddet}(A):=\operatorname{det}(AA^*)=\operatorname{det}(A^*A).
\]
Using Lemma  \ref{HerPGdetmainprothle}, we obtain
\begin{equation}\label{DetAAeqddetAAAAeq}
  \mathrm{Det}(AA^*) = \mathrm{det}(AA^*) = \mathrm{det}(A^*A) = \mathrm{Det}(A^*A).
\end{equation}
Thus,  analogously to the definition of \(\operatorname{ddet}(\cdot)\), we introduce the double determinant for an arbitrary quaternion matrix based on the \(\mathrm{P}\)-generalized determinant.

\begin{definition}\label{doubledetPgdetdef}
For \(A\in\mathbb{H}^{n\times n}\), the double determinant associated with the \(\mathrm{P}\)-generalized determinant is defined by
\[
\operatorname{dDet}(A):=\operatorname{Det}(AA^*)=\operatorname{Det}(A^*A).
\]
\end{definition}

\begin{remark}
The double determinant in Definition~\ref{doubledetPgdetdef} is defined only for square quaternion matrices.
Although both \(\operatorname{Det}(AA^*)\) and \(\operatorname{Det}(A^*A)\) are well-defined for a rectangular matrix \(A\in\mathbb{H}^{m\times n}\) with \(m\neq n\),
they are generally not equal.
In particular, if \(\operatorname{rank}(A)=n<m\), then Corollary~\ref{rankDdrowcolrcrankTH} and \cite[Theorem 7.3]{kyrchei2012theory} give
\({\rm Det}(A^*A) \neq 0\) and \({\rm Det}(AA^*) = 0\).
\end{remark}

\begin{corollary}\label{dDETbaochichengxcor}
\begin{enumerate}[(1)]
  \item\label{dDETbaochichengxcoritem01} If \(A\in \mathbb{H}^{n \times n}\), then the following statements are equivalent:
  \begin{enumerate}[(i)]
    \item\label{dDETbaochichengxcoritem01i} \(A\) is invertible;
    \item\label{dDETbaochichengxcoritem01ii} \({\rm dDet}(A) \neq 0\);
    \item\label{dDETbaochichengxcoritem01iii} \({\rm rank}(A) = n\).
  \end{enumerate}
  \item\label{dDETbaochichengxcoritem02}
  If \(A,B\in \mathbb{H}^{n \times n}\), then
$
  {\rm dDet}(AB) = {\rm dDet}(A){\rm dDet}(B).
$
\end{enumerate}
\end{corollary}

\begin{proof}
By \cite[Theorem 8.1]{kyrchei2008cramer} and \eqref{DetAAeqddetAAAAeq}, the statements \eqref{dDETbaochichengxcoritem01i} and \eqref{dDETbaochichengxcoritem01ii} are equivalent.
Moreover, \cite[Theorem 7.3]{kyrchei2012theory} and Corollary~\ref{rankDdrowcolrcrankTH} show that \eqref{dDETbaochichengxcoritem01ii} and \eqref{dDETbaochichengxcoritem01iii} are equivalent. Hence the three statements in Item~\eqref{dDETbaochichengxcoritem01} are equivalent. Finally, combining \eqref{DetAAeqddetAAAAeq} with \cite[Theorem 7.5]{kyrchei2008cramer} yields Item~\eqref{dDETbaochichengxcoritem02}.
\end{proof}

\begin{remark}
Corollary~\ref{dDETbaochichengxcor},  Theorem~\ref{DetBgBAgAvarTh}\eqref{DetBgBAgAvarThitem1}, and Theorem~\ref{DetBBgAAgvarTh}\eqref{DetBBgAAgvarThitem1} show that
the double determinant \(\operatorname{dDet}(\cdot)\)
 satisfies all three axioms for a quaternion determinant functional \(\operatorname{d}\) \cite{aslaksen1996quaternionic,cohen2000quaternionic,dyson1972quaternion}:
\begin{description}
  \item[Axiom 1] \(\operatorname{d}(A)=0\) if and only if the matrix \(A\in\mathbb{H}^{n\times n}\) is singular;
  \item[Axiom 2] \(\operatorname{d}(AB)=\operatorname{d}(A)\operatorname{d}(B)\) for all \(A,B\in\mathbb{H}^{n\times n}\);
  \item[Axiom 3] \(\operatorname{d}(P_{ij}(b)A)=\operatorname{d}(AP_{ij}(b))=\operatorname{d}(A)\)
     for all \(A\in\mathbb{H}^{n\times n}\) and \(b\in\mathbb{H}\).
\end{description}
\end{remark}

Using the double determinant \(\operatorname{dDet}(\cdot)\), we derive
 component-wise representations  for the inverse of any invertible quaternion matrix in terms of \(\mathrm{P}\)-generalized determinants.

\begin{theorem}\label{dDetaquainvedecreTH}
Let $A \in \mathbb{H}^{n \times n}$ with ${\rm dDet}(A) \neq 0$.
For all \(i,j\in[n]\),
let \(\mathrm{P}_{l\text{-}i}\) and \(\mathrm{P}_{r\text{-}j}\) be the given rules on \(S_n\), and define
\begin{equation*}
            {\mathfrak{R}}_A= \frac{1}{{\rm dDet}(A)}\left[
            \begin{array}{cccc}
          {\mathfrak{R}}_{11}  & {\mathfrak{R}}_{21}  & \ldots  & {\mathfrak{R}}_{n1} \\
          {\mathfrak{R}}_{12} & {\mathfrak{R}}_{22}  & \ldots  & {\mathfrak{R}}_{n2}\\
       \vdots  & \vdots  &   \ldots  & \vdots\\
        {\mathfrak{R}}_{1n}  & {\mathfrak{R}}_{2n}  & \ldots  & {\mathfrak{R}}_{nn}
            \end{array}
          \right]           \quad \text{and} \quad   {\mathfrak{L}}_A   =\frac{1}{{\rm dDet}(A)}
           \left[
            \begin{array}{cccc}
          {\mathfrak{L}}_{11}  & {\mathfrak{L}}_{21}  & \ldots  & {\mathfrak{L}}_{n1} \\
          {\mathfrak{L}}_{12} & {\mathfrak{L}}_{22}  & \ldots  & {\mathfrak{L}}_{n2}\\
       \vdots  & \vdots  &   \ldots  & \vdots\\
        {\mathfrak{L}}_{1n}  & {\mathfrak{L}}_{2n}  & \ldots  & {\mathfrak{L}}_{nn}
            \end{array}
          \right]   ,
\end{equation*}
where
\begin{equation*}
\mathfrak{R} _{ij} ={\rm GDet}_{{\rm P}_{l\text{-}i}} \big( (AA^*)_{i \cdot}(a_{j \cdot}^*) \big)  \quad \text{and} \quad  \mathfrak{L}_{ij} = {\rm GDet}_{{\rm P}_{r\text{-}j}} \big(  (A^*A)_{\cdot j}(a_{\cdot i}^*)\big).
\end{equation*}
Then
\begin{equation*}
  A^{-1} =     {\mathfrak{L}}_A =    {\mathfrak{R}}_A .
\end{equation*}
Moreover,
\begin{equation}\label{dDetaquainvedecreTHmaineq02}
  {\rm dDet}(A) =
   \sum \limits_{t=1}^{n} a_{kt} {\mathfrak{R}}_{kt}= \sum \limits_{t=1}^{n}  {\mathfrak{L}}_{tk}a_{tk},
\end{equation}
where $k\in [n]$.
\end{theorem}

\begin{proof}
For any $i\in [n]$, let $\operatorname{P}_{l\text{-}i}=\operatorname{P}_{l\text{-}i}(S_n, \operatorname{first}, \prec_{\triangleleft})$, and for each $j \in  [n]_i:=[n]\setminus\{i\}$,
let
\begin{equation}\label{Pliplsnifsxs2eq}
  \operatorname{P}_{l\text{-}j}^{(i)}=
\operatorname{P}_{l\text{-}j}(S_{n_1},\operatorname{first}_2,\prec_{\triangleleft_2})
\end{equation}
 be a left-order cycle representation rule, where
  $S_{n_1}$ is the symmetric group on $[n]_i$,
\[
\operatorname{first}_2(Q)=
\begin{cases}
j, & j\in\operatorname{set}(Q),\\
\operatorname{first}(Q), & j\notin\operatorname{set}(Q),
\end{cases}
\]
for all $Q\in\mathcal{Q}([n]_i)$,
and $\triangleleft_2$ is a strict total order on $[n]_i$ defined by
\[
j\triangleleft_2 d_1\triangleleft_2 d_2\triangleleft_2\cdots\triangleleft_2 d_{n-2},
\]
for  $d_1,d_2,\ldots,d_{n-2}\in[n]_i \setminus \{j\}$ satisfying
$d_1\triangleleft d_2\triangleleft\cdots\triangleleft d_{n-2}$.
Then, since ${\rm dDet}(A) \neq 0$, using \cite[Theorem 5.1]{PgeneralizedDetQuat}, we obtain
 \begin{equation}\label{AAginRAAgeq}
    ({AA^*})^{-1}= \alpha\left[
            \begin{array}{cccc}
          R_{11}  & R_{21}  & \ldots  & R_{n1} \\
          R_{12} & R_{22}  & \ldots  & R_{n2}\\
       \vdots  & \vdots  &   \ldots  & \vdots\\
        R_{1n}  & R_{2n}  & \ldots  & R_{nn}
            \end{array}
          \right],
\end{equation}
where $\alpha= \frac{1}{{\rm Det}({AA^*})}$,
\begin{equation*}
   R_{ij}=\begin{cases}
   {\rm Det}({(AA^*)}^{i,i}), & j=i,\\[1mm]
  -{\rm GDet}_{{\rm P}_{l\text{-}j}^{(i)}}\left(({(AA^*)}_{\cdot j}(b_{\cdot i}))^{i,i}\right), & j\neq i,
    \end{cases}
\end{equation*}
and $b_{\cdot i} \in \mathbb{H}^{n \times 1}$ is the $i$-th column vector of $AA^*$ for $i=1,\ldots,n$.
Let $A=[a_{xy}]$ for all $x,y \in [n]$.
Using \cite[Theorem 3.6]{PgeneralizedDetQuat}, we obtain
\begin{equation}\label{GdetAAgeqgatjRiteq}
 {\rm GDet}_{{\rm P}_{l\text{-}i}}  (AA^*)_{i \cdot}(a_{j \cdot}^*) = \sum \limits_{t=1}^{n} \overline{a_{tj}} R_{it}.
\end{equation}
Since \(A\) is invertible by Item~\eqref{dDETbaochichengxcoritem01} of Corollary~\ref{dDETbaochichengxcor},
combining \eqref{AAginRAAgeq} with \eqref{GdetAAgeqgatjRiteq} yields
\begin{align*}
A^{-1}&= A^*(AA^*)^{-1} \\
          &=
\alpha
  \left[
            \begin{array}{ccc}
          \overline{a_{11}}    & \ldots  & \overline{a_{n1}} \\
       \vdots  &   \ldots  & \vdots\\
        \overline{a_{1n}}   & \ldots  & \overline{a_{nn}}
            \end{array}
          \right]
 \left[
            \begin{array}{ccc}
          R_{11}   & \ldots  & R_{n1} \\
       \vdots  &   \ldots  & \vdots\\
        R_{1n}   & \ldots  & R_{nn}
            \end{array}
          \right]
          =
       \alpha
             \left[
            \begin{array}{ccc}
 \sum \limits_{t=1}^{n} \overline{a_{t1}} R_{1t}  &
    \ldots &  \sum \limits_{t=1}^{n} \overline{a_{t1}} R_{nt}\\
        \vdots  &   \ldots  & \vdots\\
    \sum \limits_{t=1}^{n} \overline{a_{tn}} R_{1t}  & \ldots &  \sum \limits_{t=1}^{n} \overline{a_{tn}} R_{nt}
            \end{array}
          \right]
  \\
          &=\alpha
             \left[
            \begin{array}{ccc}
 {\rm GDet}_{{\rm P}_{l\text{-}1}}  (AA^*)_{1 \cdot}(a_{1 \cdot}^*) &  \ldots&  {\rm GDet}_{{\rm P}_{l\text{-}n}}  (AA^*)_{n \cdot}(a_{1 \cdot}^*) \\
         \vdots  &   \ldots  & \vdots\\
       {\rm GDet}_{{\rm P}_{l\text{-}1}}  (AA^*)_{1 \cdot}(a_{n \cdot}^*) &  \ldots&  {\rm GDet}_{{\rm P}_{l\text{-}n}}  (AA^*)_{n \cdot}(a_{n \cdot}^*)
            \end{array}
          \right]
          =\alpha
\left[
            \begin{array}{ccc}
          {\mathfrak{R}}_{11}  &  \ldots  & {\mathfrak{R}}_{n1} \\
       \vdots  &   \ldots & \vdots\\
        {\mathfrak{R}}_{1n}  & \ldots  & {\mathfrak{R}}_{nn}
            \end{array}
          \right],
\end{align*}
which implies $A^{-1} = {\mathfrak{R}}_A$.
Furthermore,
\begin{align*}
  I_{n}& = A   {\mathfrak{R}}_A=
\alpha
\left[
            \begin{array}{ccc}
          a_{11}  &  \ldots  & a_{1n} \\
       \vdots  &   \ldots & \vdots\\
        a_{n1}  & \ldots  & a_{nn}
            \end{array}
          \right]
  \left[
            \begin{array}{ccc}
          {\mathfrak{R}}_{11}  &  \ldots  & {\mathfrak{R}}_{n1} \\
       \vdots  &   \ldots & \vdots\\
        {\mathfrak{R}}_{1n}  & \ldots  & {\mathfrak{R}}_{nn}
            \end{array}
          \right]
    =    \alpha  \left[
            \begin{array}{ccc}
       \sum\limits_{t=1}^{n} a_{1t}  {\mathfrak{R}}_{1t}  &  \ldots  &         \sum\limits_{t=1}^{n} a_{1t}  {\mathfrak{R}}_{nt} \\
       \vdots  &   \ldots &       \vdots  \\
        \sum\limits_{t=1}^{n} a_{nt}  {\mathfrak{R}}_{1t}  &  \ldots  &      \sum\limits_{t=1}^{n} a_{nt}  {\mathfrak{R}}_{nt} \\
            \end{array}
          \right],
\end{align*}
which shows that for any $k\in [n]$, we have
\begin{equation*}
   \sum \limits_{t=1}^{n} a_{kt} {\mathfrak{R}}_{kt} =   {\rm Det}({AA^*}) =  {\rm dDet}(A).
\end{equation*}
Similarly, one can prove that $A^{-1} ={\mathfrak{L}}_A$ and
$  {\rm dDet}(A) =
\sum \limits_{t=1}^{n}  {\mathfrak{L}}_{tk}a_{tk}$
for any $k\in [n]$.
\end{proof}

\begin{remark}
In particular, if the rules \(\mathrm{P}_{l\text{-}i}\) and \(\mathrm{P}_{r\text{-}j}\) satisfy Items~\eqref{PGdetspecaseremarkitem4} and~\eqref{PGdetspecaseremarkitem5} of Remark~\ref{varimdetremark}, respectively, then Theorem~\ref{dDetaquainvedecreTH} directly reduces to \cite[Theorem 8.1]{kyrchei2008cramer}.
\end{remark}

\begin{remark}
Under the hypotheses of Theorem~\ref{dDetaquainvedecreTH},
and following \cite[Remark 11.1]{kyrchei2012theory} and \cite[Remark 5.2]{PgeneralizedDetQuat},
we call the matrix
\[
\operatorname{dAdj}_{\mathrm{P}}(A):=
\begin{bmatrix}
\mathfrak{R}_{11} & \mathfrak{R}_{21} & \cdots & \mathfrak{R}_{n1}\\
\mathfrak{R}_{12} & \mathfrak{R}_{22} & \cdots & \mathfrak{R}_{n2}\\
\vdots & \vdots & \ddots & \vdots\\
\mathfrak{R}_{1n} & \mathfrak{R}_{2n} & \cdots & \mathfrak{R}_{nn}
\end{bmatrix}
=
\begin{bmatrix}
\mathfrak{L}_{11} & \mathfrak{L}_{21} & \cdots & \mathfrak{L}_{n1}\\
\mathfrak{L}_{12} & \mathfrak{L}_{22} & \cdots & \mathfrak{L}_{n2}\\
\vdots & \vdots & \ddots & \vdots\\
\mathfrak{L}_{1n} & \mathfrak{L}_{2n} & \cdots & \mathfrak{L}_{nn}
\end{bmatrix}
\]
the double \(\mathrm{P}\)-adjoint matrix of \(A\). Furthermore,
\[
\operatorname{dAdj}_{\mathrm{P}}(A)A
=A\operatorname{dAdj}_{\mathrm{P}}(A)
=\operatorname{dDet}(A)I_n.
\]
\end{remark}

Building on the above discussion of the double determinant,
the following corollary relates the double determinant to several classical quaternion determinants, namely the Moore, Xie, Chen, Study, and Dieudonn\'{e} determinants, as well as to quasideterminants. The latter were introduced by Gelfand and Retakh \cite{gelfand1991determinants,gelfand1992theory} by associating \(n^2\) determinant quantities with an \(n\times n\) matrix over a noncommutative ring. Specifically, for \(i,j\in[n]\) and an invertible matrix \(A\in\mathbb{H}^{n\times n}\), the \((i,j)\)-th quasideterminant of \(A\), denoted by \(|A|_{ij}\), is defined by
\[
|A|_{ij}=b_{ji}^{-1},
\]
where \(b_{ji}\) is the \((j,i)\)-entry of \(A^{-1}\in\mathbb{H}^{n\times n}\).

\begin{corollary}\label{dDetdderMdetCdetSdetDdetquaTH}
Let $A\in \mathbb{H}^{n\times n}$.
\begin{enumerate}[(1)]
  \item Then
  \begin{equation*}
      \begin{aligned}
   & {\rm dDet}(A) =  {\rm ddet}(A)= {\rm Mdet}(A^*A)= {\rm Xdet}(A^*A)= {\rm Cdet}(A^*A)={\rm Sdet}(A)=({\rm Ddet}(A))^2\\
  = &
     {\rm dDet}(A^*)={\rm ddet}(A^*)={\rm Mdet}(AA^*)= {\rm Xdet}(AA^*)= {\rm Cdet}(AA^*)={\rm Sdet}(A^*)=({\rm Ddet}(A^*))^2.
  \end{aligned}
  \end{equation*}

  \item\label{dDetdderMdetCdetSdetDdetquaTHitem02}
If \(A\) is invertible, then, for \(i,j\in[n]\), the quasideterminant \(|A|_{ij}\) of \(A\) is given by
  \begin{equation*}
    |A|_{ij}  =
       {\rm dDet}(A)  \,\frac{\overline{{\rm Det}_{{\rm P}_{l\text{-}i}}  \big( (AA^*)_{i \cdot}(a_{j \cdot}^*) \big) }}{|{\rm Det}_{{\rm P}_{l\text{-}i}} \big( (AA^*)_{i \cdot}(a^*_{j\cdot})  \big)|^2}=
 {\rm dDet}(A)  \, \frac{\overline{{\rm Det}_{{\rm P}_{r\text{-}j}} \big( (A^*A)_{\cdot j }(a^*_{\cdot i}) \big)}}{|{\rm Det}_{{\rm P}_{r\text{-}j}} \big( (A^*A)_{\cdot j }(a^*_{\cdot i}) \big)|^2 }.
  \end{equation*}

\end{enumerate}
\end{corollary}

\begin{proof}
\begin{enumerate}[(1)]
\item
This result follows immediately from \eqref{DetAAeqddetAAAAeq}, Lemma~\ref{HerPGdetmainprothle}, and  \cite[Remark 7.2]{kyrchei2008cramer}.

  \item
  Since \(A\) is invertible, writing \(A^{-1}=[b_{ij}]\in\mathbb{H}^{n\times n}\) and applying Theorem~\ref{dDetaquainvedecreTH}, we have
        \begin{equation*}
          b_{ij}=
           \frac{{\rm GDet}_{{\rm P}_{l\text{-}j}} \big( (AA^*)_{j \cdot}(a_{i \cdot}^*) \big) }{{\rm dDet}(A)}  =  \frac{{\rm GDet}_{{\rm P}_{r\text{-}i}} \big(  (A^*A)_{\cdot i}(a_{\cdot j}^*)\big)}{{\rm dDet}(A)} .
        \end{equation*}
        Thus, since \(\operatorname{dDet}(A)\in\mathbb{R}\), the \((i,j)\)-th quasideterminant of \(A\) is
        \begin{align*}
   |A|_{ij} = b_{ji}^{-1} &=
\left(\frac{{\rm GDet}_{{\rm P}_{l\text{-}i}} \big( (AA^*)_{i \cdot}(a_{j \cdot}^*) \big) }{{\rm dDet}(A)}\right)^{-1}= {\rm dDet}(A)  \,\frac{\overline{{\rm Det}_{{\rm P}_{l\text{-}i}}  \big( (AA^*)_{i \cdot}(a_{j \cdot}^*) \big) }}{|{\rm Det}_{{\rm P}_{l\text{-}i}} \big( (AA^*)_{i \cdot}(a^*_{j\cdot})  \big)|^2}\\
&= \left(   \frac{{\rm GDet}_{{\rm P}_{r\text{-}j}} \big(  (A^*A)_{\cdot j}(a_{\cdot i}^*)\big)}{{\rm dDet}(A)}  \right)^{-1}  =
 {\rm dDet}(A)  \, \frac{\overline{{\rm Det}_{{\rm P}_{r\text{-}j}} \big( (A^*A)_{\cdot j }(a^*_{\cdot i}) \big)}}{|{\rm Det}_{{\rm P}_{r\text{-}j}} \big( (A^*A)_{\cdot j }(a^*_{\cdot i}) \big)|^2 }.
        \end{align*}
\end{enumerate}
\end{proof}

\begin{remark}
In particular, if the rules \(\mathrm{P}_{l\text{-}i}\) and \(\mathrm{P}_{r\text{-}j}\) satisfy Items~\eqref{PGdetspecaseremarkitem4} and~\eqref{PGdetspecaseremarkitem5} of Remark~\ref{varimdetremark}, respectively, then Item~\eqref{dDetdderMdetCdetSdetDdetquaTHitem02} of Corollary~\ref{dDetdderMdetCdetSdetDdetquaTH}  yields \cite[Theorem 12.1]{kyrchei2012theory} directly.
\end{remark}

\section{Cramer's rule for quaternion equations}\label{APPsection}
In this section, we use the double determinant to derive Cramer's rule for the solutions to several classical linear  equations and matrix equations over quaternions.

\begin{theorem}\label{AxeqblineareqTH}
Let $A \in \mathbb{H}^{n\times n}$ with ${\rm dDet}(A) \neq 0$, and let $b \in \mathbb{H}^{n\times 1}$. Then the unique solution $x=[x_1,x_2,\ldots,x_n]^T \in \mathbb{H}^{n\times 1} $ of the linear equation
\begin{equation*}
  Ax=b
\end{equation*}
has components
\begin{equation*}
  x_j= \frac{{\rm GDet}_{{\rm P}_{r\text{-}j}}(A^*A)_{\cdot j}(\beta)}{{\rm dDet}(A)},
\end{equation*}
where $\beta=A^*b$ and $j=1,2,\ldots,n$.
\end{theorem}

\begin{proof}
Let $b= [b_1,b_2,\ldots,b_n]^{T}$ and
\begin{equation*}
  {\rm dAdj}_{\rm P}(A)=
  \left[
            \begin{array}{cccc}
          {\mathfrak{L}}_{11}  & {\mathfrak{L}}_{21}  & \ldots  & {\mathfrak{L}}_{n1} \\
          {\mathfrak{L}}_{12} & {\mathfrak{L}}_{22}  & \ldots  & {\mathfrak{L}}_{n2}\\
       \vdots  & \vdots  &   \ldots  & \vdots\\
        {\mathfrak{L}}_{1n}  & {\mathfrak{L}}_{2n}  & \ldots  & {\mathfrak{L}}_{nn}
            \end{array}
          \right],
\end{equation*}
where $\mathfrak{L}_{ij} = {\rm GDet}_{{\rm P}_{r\text{-}j}} \big(  (A^*A)_{\cdot j}(a_{\cdot i}^*)\big)$ with $i,j\in[n]$.
For $j\in[n]$, using Lemma \ref{basicproPgdecThle},   we obtain
\begin{align*}
\sum_{k=1}^n \mathfrak{L}_{kj}b_{k}&=
\sum_{k=1}^n
{\rm GDet}_{{\rm P}_{r\text{-}j}}
\big(  (A^*A)_{\cdot j}(a_{\cdot k}^*)\big)  \, b_{k}=\sum_{k=1}^n
{\rm GDet}_{{\rm P}_{r\text{-}j}}
 \big(  (A^*A)_{\cdot j}(a_{\cdot k}^* \cdot b_{k}) \big)
\\&=  {\rm GDet}_{{\rm P}_{r\text{-}j}} (A^*A)_{\cdot j}\left(  \sum_{k=1}^n a_{\cdot k}^* \cdot b_{k}\right)
=  {\rm GDet}_{{\rm P}_{r\text{-}j}} (A^*A)_{\cdot j}( A^*b)
=  {\rm GDet}_{{\rm P}_{r\text{-}j}} (A^*A)_{\cdot j}( \beta ).
\end{align*}
Therefore, using Theorem \ref{dDetaquainvedecreTH}, we have
\begin{align*}
x&=A^{-1}b = \frac{1}{{\rm dDet}(A)} \cdot {\rm dAdj}_{\rm P}(A) \cdot b
=
\left[
  \begin{array}{ccc}
 \frac{{\rm GDet}_{{\rm P}_{r\text{-}1}} (A^*A)_{\cdot 1}( \beta )}{{\rm dDet}(A)} &  \ldots &
 \frac{{\rm GDet}_{{\rm P}_{r\text{-}n}} (A^*A)_{\cdot n}( \beta )}{{\rm dDet}(A)}
  \end{array}
\right]^T,
\end{align*}
which completes the proof.
\end{proof}

\begin{theorem}\label{xAeqblineareqTH}
Let $A \in \mathbb{H}^{n\times n}$ with ${\rm dDet}(A) \neq 0$, and let $b \in \mathbb{H}^{1\times n}$. Then the unique solution $x=[x_1,x_2,\ldots,x_n] \in \mathbb{H}^{1\times n}$ of the linear equation
\begin{equation*}
  xA=b
\end{equation*}
has components
\begin{equation*}
  x_i= \frac{{\rm GDet}_{{\rm P}_{l\text{-}i}}(AA^*)_{i\cdot}(\gamma)}{{\rm dDet}(A)},
\end{equation*}
where $\gamma=bA^*$ and $i=1,2,\ldots,n$.
\end{theorem}

\begin{proof}
The proof is similar to that of Theorem \ref{AxeqblineareqTH} and is therefore omitted.
\end{proof}

\begin{theorem}\label{AXeqBmatrixeqTH}
Let $A \in \mathbb{H}^{m\times m}$ with ${\rm dDet}(A) \neq 0$, and let $B \in \mathbb{H}^{m\times n}$. Then the unique solution $X=[x_{ij}] \in \mathbb{H}^{m \times n}$ of the matrix equation
\begin{equation*}
  AX=B
\end{equation*}
has components
\begin{equation*}
  x_{ij}= \frac{{\rm GDet}_{{\rm P}_{r\text{-}i}}(A^*A)_{\cdot i}(\beta_{j})}{{\rm dDet}(A)},
\end{equation*}
where $\beta_{j}=A^*B_{\cdot j}$, $i=1,\ldots,m$, and $j=1,\ldots,n$.
\end{theorem}

\begin{proof}
The matrix equation \(AX=B\) is equivalent to the linear  equations
\[
AX_{\cdot j}=B_{\cdot j}, \quad j=1,\dots,n.
\]
Thus, applying Theorem~\ref{AxeqblineareqTH} to each  linear equation  yields the result.
\end{proof}

\begin{theorem}\label{XAeqBmatrixeqTH}
Let $A \in \mathbb{H}^{n\times n}$ with ${\rm dDet}(A) \neq 0$, and let $B \in \mathbb{H}^{m\times n}$. Then the unique solution $X=[x_{ij}] \in \mathbb{H}^{m \times n}$ of the matrix equation
\begin{equation*}
  XA=B
\end{equation*}
has components
\begin{equation*}
  x_{ij}= \frac{{\rm GDet}_{{\rm P}_{l\text{-}j}}(AA^*)_{j \cdot}(\gamma_{i})}{{\rm dDet}(A)},
\end{equation*}
where $\gamma_{i}=B_{i \cdot} A^*$, $i=1,\ldots,m$, and $j=1,\ldots,n$.
\end{theorem}

\begin{proof}
The proof is similar to that of Theorem \ref{AXeqBmatrixeqTH} and is therefore omitted.
\end{proof}

\begin{theorem}\label{AXBeqCmatrixeqTH}
Let \(A\in\mathbb{H}^{m\times m}\) with \(\operatorname{dDet}(A)\neq 0\), \(B\in\mathbb{H}^{n\times n}\) with \(\operatorname{dDet}(B)\neq 0\), and \(C\in\mathbb{H}^{m\times n}\).
 Then the unique solution $X=[x_{ij}] \in \mathbb{H}^{m \times n}$ of the matrix equation
\begin{equation*}
  AXB=C
\end{equation*}
has components
\begin{equation*}
  x_{ij}= \frac{{\rm GDet}_{{\rm P}_{l\text{-}j}}(BB^*)_{j\cdot}(\delta_i)}{{\rm dDet}(A){\rm dDet}(B)}
   = \frac{{\rm GDet}_{{\rm P}_{r\text{-}i}}(A^*A)_{\cdot i}(\eta_j)}{{\rm dDet}(A){\rm dDet}(B)},
\end{equation*}
where, for \(i=1,\dots,m\) and \(j=1,\dots,n\),
\begin{align*}
\delta_i&=
\begin{bmatrix}
\operatorname{GDet}_{\mathrm{P}_{r\text{-}i}}\bigl((A^*A)_{\cdot i}(\widetilde{C}_{\cdot 1})\bigr) &
\cdots &
\operatorname{GDet}_{\mathrm{P}_{r\text{-}i}}\bigl((A^*A)_{\cdot i}(\widetilde{C}_{\cdot n})\bigr)
\end{bmatrix},
\\
\eta_j&=
\begin{bmatrix}
\operatorname{GDet}_{\mathrm{P}_{l\text{-}j}}\bigl((BB^*)_{j\cdot}(\widetilde{C}_{1\cdot})\bigr) &
\cdots &
\operatorname{GDet}_{\mathrm{P}_{l\text{-}j}}\bigl((BB^*)_{j\cdot}(\widetilde{C}_{m\cdot})\bigr)
\end{bmatrix}^T,
\end{align*}
and \(\widetilde{C}=A^*CB^*\).
\end{theorem}

\begin{proof}
  Let $\widehat{a_{j\cdot}} =  (A^*A)_{j \cdot}$ for $j=1,\ldots,m$, and $\widehat{b_{\cdot i}} = (BB^*)_{\cdot i}$ for $i=1,\ldots,n$. Put
\begin{equation*}
L_{A^*A}=
  \left[
            \begin{array}{ccc}
          {{L}}_{11}   & \ldots  & {{L}}_{m1} \\
       \vdots  &   \ldots  & \vdots\\
        {{L}}_{1m}   & \ldots  & {{L}}_{mm}
            \end{array}
          \right] \quad {\text{and}} \quad
R_{BB^*}=
  \left[
            \begin{array}{ccc}
          {{R}}_{11}    & \ldots  & {{R}}_{n1} \\
       \vdots    &   \ldots  & \vdots\\
        {{R}}_{1n}   & \ldots  & {{R}}_{nn}
            \end{array}
          \right],
\end{equation*}
where
\begin{align*}
  L_{ij}&=\begin{cases}
   {\rm Det}((A^*A)^{j,j}), & i=j,\\[1mm]
   -{\rm GDet}_{{\rm P}_{r\text{-}i}^{(j)}}\left(((A^*A)_{i\cdot}( \widehat{a_{j\cdot}}))^{j,j}\right), & i \in [m]\setminus\{j\},
    \end{cases} \\
   R_{ij}&=\begin{cases}
   {\rm Det}((BB^*)^{i,i}), & j=i,\\[1mm]
  -{\rm GDet}_{{\rm P}_{l\text{-}j}^{(i)}}\left(((BB^*)_{\cdot j}(   \widehat{b_{\cdot i}}))^{i,i}\right), & j \in [m]\setminus\{i\},
    \end{cases}
\end{align*}
the rules $\mathrm{P}_{l\text{-}j}^{(i)}$ are defined in \eqref{Pliplsnifsxs2eq} for any $i\in[m]$ and each $j\in [m]\setminus\{i\}$, and $\mathrm{P}_{r\text{-}i}^{(j)}$ are defined similarly for  any $j\in[m]$ and each $i\in [m]\setminus\{j\}$.
Using  \cite[Theorems 5.1, 3.6]{PgeneralizedDetQuat}, we obtain
\begin{align*}
(A^*A)^{-1}\widetilde{C}&=(A^*A)^{-1}(A^*CB^*)=  \frac{1}{{\rm dDet}(A)}  \, L_{A^*A} \, (A^*CB^*)\\
&=  \frac{1}{{\rm dDet}(A)}  \,
\left[
     \begin{array}{ccc}
    {\rm GDet}_{{\rm P}_{r\text{-}1}} (A^*A)_{\cdot 1}( \widetilde{C}_{\cdot 1})    &  \ldots &  {\rm GDet}_{{\rm P}_{r\text{-}1}} (A^*A)_{\cdot 1}( \widetilde{C}_{\cdot n})  \\
 \vdots   &  \ldots &  \vdots\\
        {\rm GDet}_{{\rm P}_{r\text{-}m}} (A^*A)_{\cdot m}( \widetilde{C}_{\cdot 1})    &  \ldots &  {\rm GDet}_{{\rm P}_{r\text{-}m}} (A^*A)_{\cdot m}( \widetilde{C}_{\cdot n})  \\
     \end{array}
   \right].
\end{align*}
Therefore,
\begin{align*}
X&= (A^*A)^{-1}(A^*CB^*)(BB^*)^{-1} = (A^*A)^{-1}\widetilde{C}(BB^*)^{-1}\\
& = \alpha  \,
\left[
     \begin{array}{ccc}
    {\rm GDet}_{{\rm P}_{r\text{-}1}} (A^*A)_{\cdot 1}( \widetilde{C}_{\cdot 1})    &  \ldots &  {\rm GDet}_{{\rm P}_{r\text{-}1}} (A^*A)_{\cdot 1}( \widetilde{C}_{\cdot n})  \\
 \vdots   &  \ldots &  \vdots\\
        {\rm GDet}_{{\rm P}_{r\text{-}m}} (A^*A)_{\cdot m}( \widetilde{C}_{\cdot 1})    &  \ldots &  {\rm GDet}_{{\rm P}_{r\text{-}m}} (A^*A)_{\cdot m}( \widetilde{C}_{\cdot n})  \\
     \end{array}
   \right]
 \left[
            \begin{array}{ccc}
          {{R}}_{11}  &  \ldots  & {{R}}_{n1} \\
       \vdots  &   \ldots  & \vdots\\
        {{R}}_{1n}   & \ldots  & {{R}}_{nn}
            \end{array}
          \right]\\
&  = \alpha  \,
\left[
            \begin{array}{ccc}
{\rm GDet}_{{\rm P}_{l\text{-}1}}(BB^*)_{ 1\cdot }(\delta_1) &  \ldots  &  {\rm GDet}_{{\rm P}_{l\text{-}n}}(BB^*)_{ n\cdot }(\delta_1)\\
       \vdots  &   \ldots  & \vdots\\
{\rm GDet}_{{\rm P}_{l\text{-}1}}(BB^*)_{ 1\cdot }(\delta_m) &  \ldots  &  {\rm GDet}_{{\rm P}_{l\text{-}n}}(BB^*)_{ n\cdot }(\delta_m)\\
            \end{array}
          \right],
\end{align*}
where $\alpha =\frac{1}{{\rm dDet}(A) {\rm dDet}(B)} $.
Similarly, one can prove that
  \begin{equation*}
  x_{ij}=
   \frac{{\rm GDet}_{{\rm P}_{r\text{-}i}}(A^*A)_{ \cdot i }(\eta_j)}{{\rm dDet}(A){\rm dDet}(B)},
\end{equation*}
where  $i=1,\ldots,m$ and $j=1,\ldots,n$.
\end{proof}

\begin{remark}
In particular, if the rules \(\mathrm{P}_{l\text{-}i}\) and \(\mathrm{P}_{r\text{-}j}\) satisfy Items~\eqref{PGdetspecaseremarkitem4} and~\eqref{PGdetspecaseremarkitem5} of Remark~\ref{varimdetremark}, respectively, then Theorems~\ref{AxeqblineareqTH}, \ref{xAeqblineareqTH}, \ref{AXeqBmatrixeqTH}, \ref{XAeqBmatrixeqTH}, and \ref{AXBeqCmatrixeqTH} reduce to \cite[Theorems 2.11, 2.12, 3.1, 3.2, 3.3]{kyrchei2010cramersq}, respectively.
\end{remark}

\section{Application to color image encryption and decryption}\label{numexpsection}
Traditional methods process the red, green, and blue (RGB) channels separately, thereby ignoring their inherent color correlation. In contrast, the quaternion representation encodes RGB as a single entity, preserving channel coupling and avoiding color distortion \cite{grigoryan2018quaternion}.
In recent years, color image encryption and decryption methods that rely on solving matrix equations have attracted increasing attention~\cite{xie2025generalized,gao2026novel}.
In this section, we design a color image encryption--decryption scheme based on Cramer's rule for the unique solution to the quaternion matrix equation \(AXB=C\) established in Theorem~\ref{AXBeqCmatrixeqTH}.

However, directly solving the quaternion matrix equation for the entire image requires a computational cost of at least \(O(mn(m!+n!))\).
Therefore, we adopt a block strategy: the whole image is divided into small blocks, and encryption and decryption are performed block by block, thereby significantly reducing the computational complexity while ensuring security.
In the block strategy, when the image dimensions are not divisible by block sizes, zero padding is applied along the right and bottom boundaries, which is subsequently removed during decryption.
To further enhance security, Arnold scrambling \cite{arnold1968ergodic} is incorporated into the encryption, with the number of iterations also included as a secret key.
Since Arnold transformation is only defined on square matrices, the original rectangular image is first zero-padded to a square before scrambling.

\begin{figure}[!htb]
\centering
\includegraphics[width=0.95\textwidth]{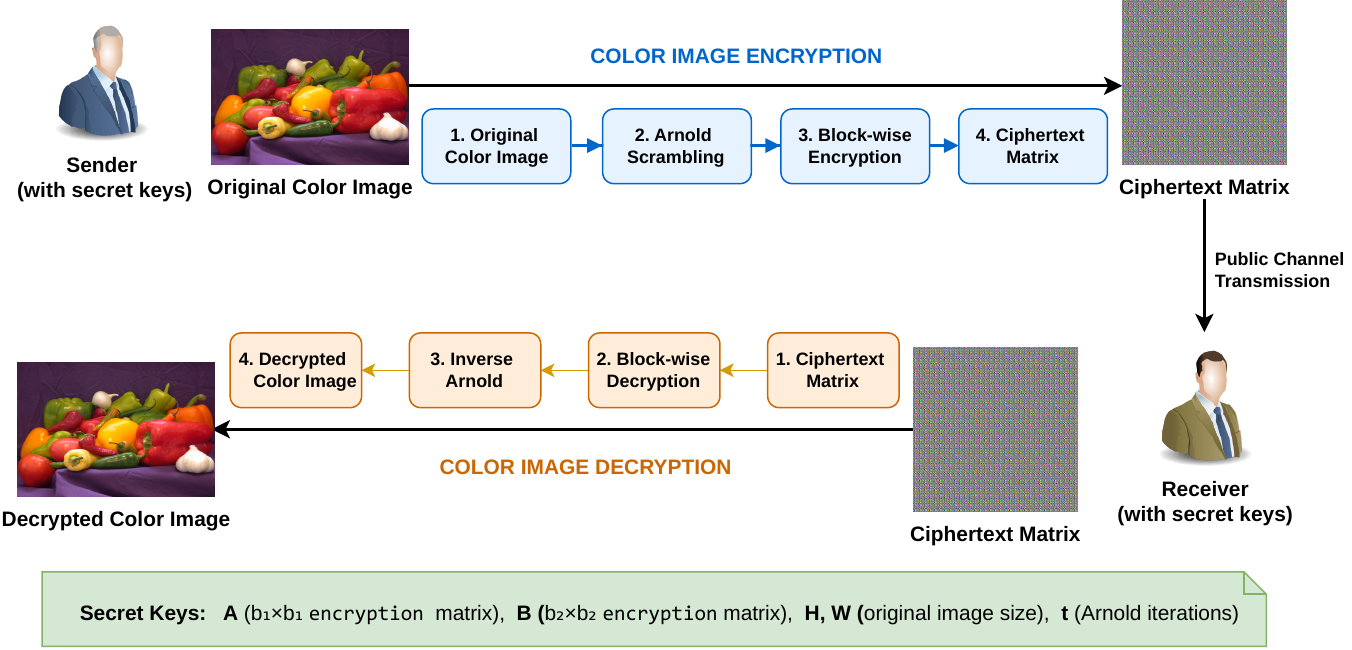}
\caption{The framework of the proposed encryption--decryption scheme}
\label{fig:framework}
\end{figure}

The encryption and decryption procedures are described in
Algorithms~\ref{alg:enc} and~\ref{alg:dec}, respectively.
The transmission protocol between the sender and receiver is detailed as follows, and the overall framework of the proposed  encryption--decryption scheme is illustrated in Figure~\ref{fig:framework}.

\begin{algorithm}[!htb]
\caption{Color image encryption}
\label{alg:enc}
\begin{algorithmic}[1]
\Require Original color image $\mathcal{X}$ of size $H\times W$;
 encryption matrices $A\in\mathbb{H}^{b_1\times b_1}$ and $B\in\mathbb{H}^{b_2\times b_2}$ with ${\rm dDet}(A)\neq 0$ and ${\rm dDet}(B)\neq 0$;
 Arnold iteration number $t$.
\Ensure Ciphertext matrix $C\in\mathbb{H}^{S_1\times S_2}$, where $S_1=\lceil S/b_1\rceil b_1$, $S_2=\lceil S/b_2\rceil b_2$, $S=\max(H,W)$.
\State Take the RGB values of $\mathcal{X}$ as the imaginary parts $\mathbf{i},\mathbf{j},\mathbf{k}$ to form the pure quaternion matrix $Y_0\in\mathbb{H}^{H\times W}$.
\State Pad $Y_0$ with zeros on the bottom and right to obtain $Y_1\in\mathbb{H}^{S\times S}$, where $S=\max(H,W)$.
\State Apply Arnold scrambling $t$ times on $Y_1$ to obtain the scrambled matrix $Y_2\in\mathbb{H}^{S\times S}$.
\State Pad $Y_2$ with zeros on the bottom and right to obtain $Y_3\in\mathbb{H}^{S_1\times S_2}$, so that $b_1\mid S_1$ and $b_2\mid S_2$.
\State Partition $Y_3$ into $L=\frac{S_1}{b_1}\cdot\frac{S_2}{b_2}$ non-overlapping $b_1\times b_2$ subblocks $\{X_\ell\}_{\ell=1}^L$.
\State Compute $C_\ell=AX_\ell B$ for each subblock $X_\ell$.
\State Concatenate all $C_\ell$ to obtain the ciphertext matrix $C\in\mathbb{H}^{S_1\times S_2}$.
\end{algorithmic}
\end{algorithm}

\begin{algorithm}[!htb]
\caption{Color image decryption}
\label{alg:dec}
\begin{algorithmic}[1]
\Require Ciphertext matrix $C\in\mathbb{H}^{S_1\times S_2}$;
\textbf{Secret keys:} encryption matrices $A\in\mathbb{H}^{b_1\times b_1}$, $B\in\mathbb{H}^{b_2\times b_2}$, original image size $H\times W$, Arnold iteration number $t$.
\Ensure Decrypted color image $\widehat{\mathcal{X}}$ of size $H\times W$.
\State Partition $C$ into $L=\frac{S_1}{b_1}\cdot\frac{S_2}{b_2}$ non-overlapping $b_1\times b_2$ subblocks $\{C_\ell\}_{\ell=1}^L$.
\State For each subblock $C_\ell$, solve $A\widehat{X}_\ell B=C_\ell$ by Cramer's rule in Theorem~\ref{AXBeqCmatrixeqTH}.
\State Concatenate all $\widehat{X}_\ell$ to obtain the padded quaternion matrix $Y_3\in\mathbb{H}^{S_1\times S_2}$.
\State Crop the bottom $S_1-S$ rows and right $S_2-S$ columns of $Y_3$ to recover the square matrix $Y_2\in\mathbb{H}^{S\times S}$, where $S=\max(H,W)$.
\State Apply inverse Arnold scrambling $t$ times on $Y_2$ to obtain the unscrambled matrix $Y_1\in\mathbb{H}^{S\times S}$.
\State Crop the bottom $S-H$ rows and right $S-W$ columns of $Y_1$ to obtain the pure quaternion matrix $Y_0\in\mathbb{H}^{H\times W}$.
\State Extract the imaginary parts $\mathbf{i},\mathbf{j},\mathbf{k}$ of $Y_0$ as the RGB channels to obtain the
decrypted color image
$\widehat{\mathcal{X}}$ of size $H\times W$.

\end{algorithmic}
\end{algorithm}

\noindent\textbf{Transmission protocol:} Before transmission,
the sender and receiver pre-share the following secret keys through a secure channel: the encryption matrices $A\in\mathbb{H}^{b_1\times b_1}$ and $B\in\mathbb{H}^{b_2\times b_2}$,  the original image size  $H\times W$, and the Arnold iteration number $t$.
Since these keys are small in size, their secure distribution is low-cost. During encryption, the sender executes Algorithm~\ref{alg:enc} and transmits the data-intensive ciphertext matrix $C$ over a public channel.
Using the shared keys, the receiver then executes Algorithm~\ref{alg:dec} to obtain the decrypted color image $\widehat{\mathcal{X}}$.
Even if an adversary intercepts the entire ciphertext matrix $C$,
without the full set of secret keys they cannot recover the original image.

To quantitatively evaluate the encryption effect, the information entropy ~\cite{netravali1995digital} of the ciphertext matrix is computed by normalizing its three imaginary-part channels to the pixel range $0$--$255$.
 The closer the entropy value is to the theoretical maximum $8$ bits, the stronger the encryption effect.
For decryption quality, the mean squared error (MSE) and peak signal-to-noise ratio (PSNR) are adopted~\cite{wu2013local}.
A smaller MSE indicates more accurate decryption,
while a higher PSNR indicates better quality;
above $40$ dB the decrypted image is visually indistinguishable from the original, and it is infinity upon perfect recovery.

\begin{figure}[!ht]
\centering
\includegraphics[width=0.95\textwidth,trim=3cm 2cm 2cm 0.1cm,clip]{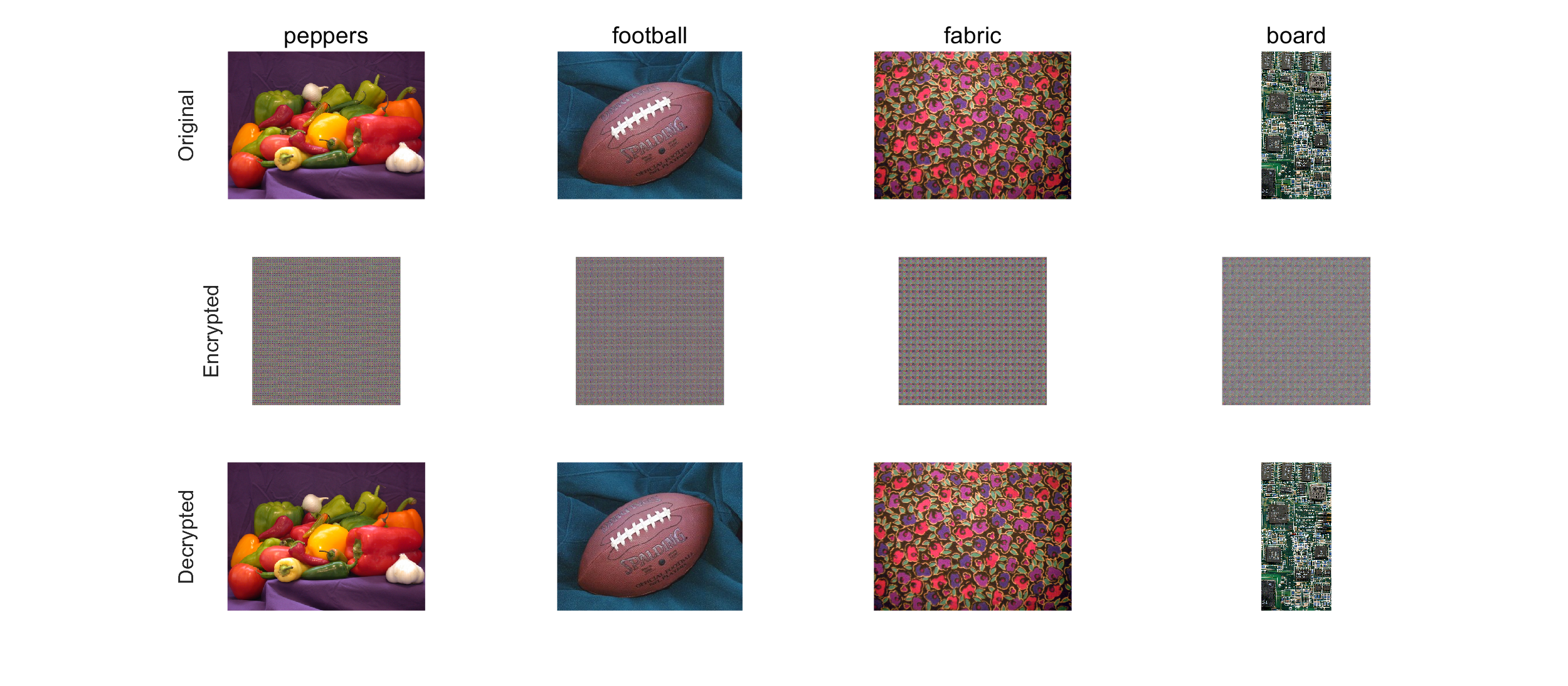}
\caption{Encryption and decryption results on four test images. }
\label{fig:results}
\end{figure}

Four test color images, namely \textit{peppers.png}, \textit{football.jpg}, \textit{fabric.png}, and \textit{board.tif}, are selected from the MATLAB Image Processing Toolbox to evaluate the proposed algorithm.
 The row-column determinant is adopted as a particular $\rm P$-generalized determinant for numerical experiments.
  The secret keys for the four test images are set as follows:
\begin{enumerate}[(1)]
  \item The encryption matrices $A \in \mathbb{H}^{b_1 \times b_1}$ and
$B \in \mathbb{H}^{b_2 \times b_2}$ are given by
\begin{align*}
A &= \left[\begin{array}{ccc}
-3 + 4\mathbf{i} - \mathbf{j} & 2 + 7\mathbf{i} - 8\mathbf{j} + 2\mathbf{k} & -9 - 7\mathbf{i} - \mathbf{j} + 9\mathbf{k} \\
9 - 10\mathbf{i} - 4\mathbf{j} + 2\mathbf{k} & -7 - 6\mathbf{i} - 4\mathbf{j} - 7\mathbf{k} & 8 - 4\mathbf{i} + 6\mathbf{j} + 10\mathbf{k} \\
5 + 10\mathbf{i} + 2\mathbf{j} - 10\mathbf{k} & -7 - 7\mathbf{i} - 3\mathbf{j} - 9\mathbf{k} & 2 + \mathbf{i} - 6\mathbf{j} + 6\mathbf{k}
\end{array}\right],
\\
B &=  \left[ \begin{array}{ccc}
-4 + 3\mathbf{i} + 2\mathbf{j} - 3\mathbf{k} & -1 + \mathbf{i} - 6\mathbf{j} - 8\mathbf{k} & -10 + 6\mathbf{i} - 2\mathbf{j} + 10\mathbf{k} \\
-8 - 4\mathbf{i} + 9\mathbf{j} - 5\mathbf{k} & -8 - 7\mathbf{i} - 10\mathbf{j} + 6\mathbf{k} & 9 + 9\mathbf{i} - 5\mathbf{j} + 6\mathbf{k} \\
4 - 9\mathbf{j} + \mathbf{k} & 10\mathbf{i} - 4\mathbf{j} - 9\mathbf{k} & -5 + 8\mathbf{i} + 7\mathbf{j} - 6\mathbf{k}
\end{array}\right],
\end{align*}
where $b_1 = 3$ and $b_2 = 3$ are the row and column block sizes, respectively.
\item  The original image sizes $(H, W)$ for the four test images are
$384 \times 512$ (peppers), $256 \times 320$ (football),
$480 \times 640$ (fabric), and $648 \times 306$ (board), respectively.

\item  The number of Arnold scrambling iterations is set to $10$.

\end{enumerate}
To visually demonstrate the encryption effect, the imaginary parts of the ciphertext matrix $C$ are displayed as an encrypted image.
All experiments are implemented in MATLAB R2019a on a laptop with an AMD Ryzen 7 5800H CPU (3.20 GHz) and 16 GB RAM.
The encryption and decryption results are shown in Figure~\ref{fig:results}, and the quantitative metrics are listed in Table~\ref{tabmetrics}.

\begin{table}[!ht]
\centering
\setlength{\tabcolsep}{20pt}
\renewcommand{\arraystretch}{1.3}
\begin{tabular}{ccccc}
\hline
\textbf{Image} & \textbf{Entropy (bits)} & \textbf{MSE} & \textbf{PSNR (dB)} \\
\hline
peppers  & 6.143 & 0.0000 & Inf \\
football & 5.869 & 0.0000 & Inf \\
fabric   & 6.348 & 0.0000 & Inf \\
board    & 5.687 & 0.0000 & Inf \\
\hline
\end{tabular}
\caption{Encryption and decryption metrics for the four test images}
\label{tabmetrics}
\end{table}

From Figure~\ref{fig:results}, the encrypted images appear as random noise, completely concealing the contours and color information of the originals; meanwhile, the decrypted images are nearly identical to the originals and visually indistinguishable.
Furthermore, Table~\ref{tabmetrics} shows that the information entropies of the encrypted images are reasonably high, indicating that the pixel distribution is substantially randomized.
The MSE values are zero and the PSNR values are infinite for all four images, demonstrating that the proposed algorithm achieves exact decryption and  almost perfectly recovers the original color images.

\section{Conclusions}\label{conclussection}
Using the $\rm P$-generalized determinant, we define the double determinant and apply it to the determinant representation of  matrix inverses, as well as to Cramer's rule for the unique solutions of
the quaternion  linear equations \(Ax=b\) and \(xA=b\), as well as  of the quaternion matrix equations \(AX=B\), \(XA=B\), and \(AXB=C\). Furthermore, based on the solution of $AXB=C$, we design
the color image encryption--decryption scheme. Numerical experiments demonstrate that the encryption achieves satisfactory results and the decryption nearly perfectly recovers the original color images.
Quaternion determinants constitute an important part of the fundamental theory of quaternion algebra. We expect that the $\rm P$-generalized determinant and the double determinant will find further applications in classical areas where quaternions have been successfully employed, such as color image processing \cite{grigoryan2018quaternion}, computer graphics\cite{goldman2010rethinking,vince2011quaternions}, aerospace engineering\cite{kuipers1999quaternions}, and quantum mechanics \cite{adler1995quaternionic}.

\section*{Funding}
This work was  supported by the National Natural Science Foundation of China [grant number 12371023].

\section*{Conflict of interest}
The authors declare no conflict of interest.

\section*{Author contributions}
\noindent\textbf{Jiale Gao:} Conceptualization, Methodology, Formal analysis, Investigation, Writing -- Original Draft, Writing -- Review \& Editing.

\noindent \textbf{Qing-Wen Wang:} Methodology, Funding acquisition, Supervision, Writing -- Review \& Editing.


\begin{thebibliography}{10}

\bibitem{adler1995quaternionic}
S.~L. Adler.
\newblock {\em Quaternionic Quantum Mechanics and Quantum Fields}.
\newblock Oxford University Press, New York, NY, USA, 1995.

\bibitem{arnold1968ergodic}
V.~I. Arnold and A.~Avez.
\newblock {\em Ergodic Problems of Classical Mechanics}.
\newblock W. A. Benjamin, New York, NY, USA, 1968.

\bibitem{aslaksen1996quaternionic}
H.~Aslaksen.
\newblock Quaternionic determinants.
\newblock {\em Math. Intelligencer}, 18(3):57--65, 1996.

\bibitem{cayley1845quat}
A.~Cayley.
\newblock {XIII}. {On} certain results relating to quaternions.
\newblock {\em Philos. Mag.}, 26(171):141--145, 1845.

\bibitem{chen1991definition}
L.~Chen.
\newblock Definition of determinant and {Cramer} solutions over the quaternion
  field.
\newblock {\em Acta Math. Sin. New Ser.}, 7(2):171--180, 1991.

\bibitem{cohen2000quaternionic}
N.~Cohen and S.~De~Leo.
\newblock The quaternionic determinant.
\newblock {\em Electron. J. Linear Algebra}, 7:100--111, 2000.

\bibitem{dieudonne1943determinants}
J.~Dieudonn{\'e}.
\newblock Les d{\'e}terminants sur un corps non commutatif.
\newblock {\em Bull. Soc. Math. Fr.}, 71:27--45, 1943.
\newblock In French.

\bibitem{dyson1972quaternion}
F.~J. Dyson.
\newblock Quaternion determinants.
\newblock {\em Helv. Phys. Acta}, 45:289--302, 1972.

\bibitem{PgeneralizedDetQuat}
J.~Gao and Q.-W. Wang.
\newblock {$\rm P$}-generalized determinant over quaternions and its
  applications.
\newblock arXiv, 2026.
\newblock arXiv:2610.41644.

\bibitem{gao2026novel}
Z.-H. Gao, Q.-W. Wang, and L.-M. Xie.
\newblock A novel {Moore--Penrose} inverse to dual quaternion matrices with
  applications.
\newblock {\em Appl. Math. Lett.}, 172:109727, 2026.

\bibitem{gelfand1991determinants}
I.~M. Gelfand and V.~S. Retakh.
\newblock Determinants of matrices over noncommutative rings.
\newblock {\em Funct. Anal. Appl.}, 25(2):91--102, 1991.

\bibitem{gelfand1992theory}
I.~M. Gelfand and V.~S. Retakh.
\newblock A theory of noncommutative determinants and characteristic functions
  of graphs.
\newblock {\em Funct. Anal. Appl.}, 26(4):231--246, 1992.

\bibitem{goldman2010rethinking}
R.~Goldman.
\newblock {\em Rethinking Quaternions: Theory and Computation}.
\newblock Synthesis Lectures on Computer Graphics and Animation, no. 13. Morgan
  \& Claypool, 2010.

\bibitem{grigoryan2018quaternion}
A.~M. Grigoryan and S.~S. Agaian.
\newblock {\em Quaternion and Octonion Color Image Processing with {MATLAB}}.
\newblock SPIE, Bellingham, WA, USA, 2018.

\bibitem{hamilton1853lectures}
W.~R. Hamilton.
\newblock {\em Lectures on Quaternions}.
\newblock Hodges and Smith, Dublin, 1853.

\bibitem{hamilton1899elements}
W.~R. Hamilton.
\newblock {\em Elements of Quaternions (Second Edition)}.
\newblock Longmans, Green, and Co., London, 1899.

\bibitem{hanson2006visualizing}
A.~J. Hanson.
\newblock {\em {Visualizing Quaternions}}.
\newblock Morgan Kaufmann, San Francisco, CA, USA, 2006.

\bibitem{hanson2024visualizing}
A.~J. Hanson.
\newblock {\em {Visualizing More Quaternions}}.
\newblock Morgan Kaufmann, Cambridge, MA, USA, 2024.

\bibitem{kuipers1999quaternions}
J.~B. Kuipers.
\newblock {\em Quaternions and Rotation Sequences: A Primer with Applications
  to Orbits, Aerospace, and Virtual Reality}.
\newblock Princeton University Press, Princeton, NJ, USA, 1999.

\bibitem{kyrchei2008cramer}
I.~I. Kyrchei.
\newblock {Cramer's rule} for quaternionic systems of linear equations.
\newblock {\em J. Math. Sci.}, 155(6):839--858, 2008.

\bibitem{kyrchei2010cramersq}
I.~I. Kyrchei.
\newblock {Cramer's} rule for some quaternion matrix equations.
\newblock {\em Appl. Math. Comput.}, 217(5):2024--2030, 2010.

\bibitem{kyrchei2012theory}
I.~I. Kyrchei.
\newblock The theory of the column and row determinants in a quaternion linear
  algebra.
\newblock In A.~R. Baswell, editor, {\em Advances in Mathematics Research},
  volume~15, pages 301--358. Nova Sci. Publ., New York, NY, USA, 2012.

\bibitem{Kyrchei2026}
I.~I. Kyrchei.
\newblock {\em Quaternion Generalized Inverses: Foundations, Theory, Problems,
  and Solutions}.
\newblock Academic Press, London, 2026.

\bibitem{moore1922determinant}
E.~H. Moore.
\newblock On the determinant of an {Hermitian} matrix of quaternionic elements.
\newblock {\em Bull. Amer. Math. Soc.}, 28:161--162, 1922.

\bibitem{netravali1995digital}
A.~N. Netravali and B.~G. Haskell.
\newblock {\em Digital Pictures: Representation, Compression, and Standards
  ({Second} Edition)}.
\newblock Springer Science+Business Media, New York, NY, USA, 1995.

\bibitem{study1920zur}
E.~Study.
\newblock Zur theorie der linearen gleichungen.
\newblock {\em Acta Math.}, 42:1--61, 1920.
\newblock In German.

\bibitem{vince2011quaternions}
J.~Vince.
\newblock {\em Quaternions for Computer Graphics (Second Edition)}.
\newblock Springer, London, 2021.

\bibitem{wu2013local}
Y.~Wu, Y.~Zhou, G.~Saveriades, S.~Agaian, J.~P. Noonan, and P.~Natarajan.
\newblock Local {Shannon} entropy measure with statistical tests for image
  randomness.
\newblock {\em Inf. Sci.}, 222:323--342, 2013.

\bibitem{xie1980centralizable}
B.~Xie.
\newblock Determinants of centralizable matrices over any skew-field.
\newblock {\em J. Jilin Univ. Natur. Sci.}, (3):1--33, 1980.
\newblock In Chinese.

\bibitem{xie1980expansion}
B.~Xie.
\newblock Expansion theorem of the determinant of self-conjugate quaternion
  matrices and its applications.
\newblock {\em Acta Math. Sin.}, 23(5):668--683, 1980.
\newblock In Chinese.

\bibitem{xie2025generalized}
L.-M. Xie, Q.-W. Wang, and Z.-H. He.
\newblock The generalized hand-eye calibration matrix equation {$AX-YB=C$} over dual quaternions.
\newblock {\em Comput. Appl. Math.}, 44:137, 2025.

\end{thebibliography}
\end{document}